\documentclass[a4paper,11pt]{article}

\usepackage{amssymb,amsthm,amsmath}

\usepackage[ruled]{algorithm}
\usepackage{algpseudocode}

\renewcommand{\algorithmicrequire}{\textbf{Input:}}
\renewcommand{\algorithmicensure}{\textbf{Output:}}

\newcommand{\VC}{V_\oC}
\newcommand{\Tr}{{\text{Tr}}}
\newcommand{\fCM}{M^{2\oC}}
\newcommand{\fCK}{K^{2\oC}}
\newcommand{\ty}{\tilde y}
\newcommand{\Span}{{\rm Span}}
\newcommand{\PR}{\oR[x]}
\newcommand{\PK}{\oK[x]}
\newcommand{\PC}{\oC[x]}

\newcommand{\LT}{{\rm LT}}
\newcommand{\oK}{\mathbb K}

\newcommand{\res}{{\text{res}}}

\newcommand{\oT}{\mathbb T}

\newcommand{\Ker}{{\rm Ker}}
\newcommand{\rank}{{\rm rank}}
\newcommand{\oR}{\mathbb R}

\newcommand{\oN}{\mathbb N}
\newcommand{\oC}{\mathbb C}
\newcommand{\oz}{\overline x}
\newcommand{\ov}{\overline v}
\newcommand{\ox}{\overline x}

\newcommand{\ZZ}{\bar \oT}

\newcommand{\oX}{\overline x}
\newcommand{\CC}{{\cal C}}

\newcommand{\XX}{{\cal X}}
\newcommand{\pBB}{\partial \BB}

\newcommand{\BB}{{\cal B}}

\newcommand{\ve}[1]{\ensuremath{{\rm vec}(#1)}}

\numberwithin{equation}{section}

\newtheorem{theorem}{\bf Theorem}[section]
 
\newtheorem{lemma}[theorem]{\bf Lemma}
 \newtheorem{claim}[theorem]{\bf Claim}
 \newtheorem{proposition}[theorem]{\bf Proposition}

 \newtheorem{corollary}[theorem]{\bf Corollary}

\theoremstyle{remark}
\newtheorem{remark}[theorem]{\bf Remark}

\newtheoremstyle{example}
  {}
  {}
  {\normalfont}
  {}
  {\bfseries}
  {}
  {9pt}
  {}
\theoremstyle{example}
\newtheorem{example}[theorem]{Example}

\begin{document}
\begin{center}{\Large \bf Semidefinite characterization and computation of
zero-dimensional real radical ideals}

\bigskip
{\Large J.-B. Lasserre\footnote{
LAAS-CNRS and Institute of Mathematics, Toulouse, France. 
{\tt lasserre@laas.fr}. Supported by the french national research agency ANR under grant NT05-3-41612}, M. Laurent\footnote{Centrum voor Wiskunde en Informatica, Kruislaan 413, 1098 SJ Amsterdam, Netherlands. {\tt M.Laurent@cwi.nl}.
Supported by the Netherlands Organization for Scientific Research grant
NWO 639.032.203 and
by ADONET, Marie Curie Research Training Network MRTN-CT-2003-504438.}
, P. Rostalski\footnote{
Automatic Control Laboratory, Physikstrasse 3, ETH Zurich, 8092 Zurich, Switzerland.
{\tt rostalski@control.ee.ethz.ch}.}}

\bigskip
\end{center}

\begin{abstract}
For an ideal $I\subseteq\oR[x]$ given by a set of generators, 
a new semidefinite characterization of its real radical $I(V_\oR(I))$ 
is presented, provided it is zero-dimensional 
(even if $I$ is not). 
Moreover we propose an algorithm using numerical linear algebra and semidefinite optimization
 techniques,
 to compute all (finitely many) points of the real variety $V_\oR(I)$ as well
 as a set of generators of the real radical ideal. The latter is obtained in the form of 
a border or  Gr\"obner basis. The algorithm is based on moment relaxations  and,
 in contrast to other existing methods, it exploits the real algebraic nature of the problem right from the beginning and avoids the computation of complex components.\\
 
\noindent
 {\bf AMS:} 14P05 13P10 12E12 12D10 90C22\\
 

\end{abstract}

\section{Introduction}	  
\label{sec::intro}

Algebraic computation over the reals is a highly relevant topic with many practical applications and, in particular, for finding {real} solutions to a system of polynomial equations.
Throughout the paper, $\PK=\oK[x_1,\ldots,x_n]$ 
denotes the ring of polynomials in $n$ variables over the field $\oK=\oR$ or $\oC$.
For an ideal $I\subseteq \PK$, 
$V_\oC(I):=\{x\in\oC^n\mid f(x)=0\ \forall f\in I\}$ and
$V_\oR(I):=\VC(I)\cap \oR^n$
denote, respectively,  the complex and real varieties of $I$ and, for $V\subseteq \oC^n$,
$I(V):=\{f\in \PK\mid f(v)=0\ \forall v\in V\}$ is the vanishing ideal of the set $V$. The ideal
$I(V_\oC(I))$ coincides with the radical ideal $\sqrt I$ of $I$ by the Nullstellensatz and
$I(V_\oR(I))$ coincides with the real radical ideal $\sqrt[\oR]{I}$ by the Real Nullstellensatz (see Section \ref{secprelideal1} for details).
The problem of finding the radical ideal $I(V_\oC(I))$ seems to be 
much better understood than that of finding the real radical ideal $I(V_\oR(I))$; see below for a brief recap on existing literature.
In this paper, we provide a new characterization of the real radical ideal $I(V_\oR(I))$ of an
ideal  $I\subseteq \oR[x]$, assuming $I$ is given by generators $h_1,\ldots,h_m \in \PR$ and $V_\oR(I)$ is finite (while $V_\oC(I)$ needs not be finite).
In addition, from this characterization, we also define a numerical algorithm based on semidefinite programming to compute the points of the (finite) variety $V_\oR(I)$ as well as a set of generators of the real radical ideal $I(V_\oR(I))$.
More generally, our results extend to the case of the so-called $S$-radical ideal $I(V_\oR(I)\cap S)$ where $S\subseteq \oR^n$ is defined by finitely many  polynomial inequalities, assuming that $V_\oR(I) \cap S$ is finite.
It turns out that a similar algorithm also works for computing $\VC(I)$ and the radical ideal $I(\VC(I))$ (assuming now $\VC(I)$ finite) 
although very good methods already exist  for this latter case.
In the remainder of the Introduction, after recalling some motivation and related literature on the problem of finding the (real) radical ideal, 
we sketch the main ingredients of our method. We already introduce some definitions but refer to 
Sections \ref{secprelideal} and \ref{secprelmom} for additional   definitions about polynomials and moment matrices.

\medskip\noindent
{\bf \em Motivation.} The main motivation of this work is to provide a characterization as well as an algorithm for finding the real variety and
the real radical of an ideal $I\subseteq \oR[x]$  that takes into account the
 specific {\it real} algebraic geometric nature of the problem. Indeed, to the best of our knowledge, most basic methods for computing 
the real variety $V_\oR(I)$ 
first compute the complex variety $\VC(I)$; for this
they  require as basic ingredients a Gr\"obner basis of $I$
and a  linear basis of the vector space $\oR[x]/I$ and thus they work under the assumption 
that $\VC(I)$ is finite.
Even if  $\VC(I)$ is finite but has many more complex elements than real ones, this may produce a large computational overhead. This is particularly important as the numbers of complex and real solutions may differ significantly as supported by the fewnomial theory of Khovanski \cite{khov}; see also the discussion in Bihan et al. \cite{bihan}, \cite{bihansottile}.
In other words, this problem of real algebraic geometry is solved via algebraic methods that do not take into account {right from the beginning} the real algebraic aspect of the problem. In contrast,
our characterization and our algorithm do {not} need knowledge of a 
Gr\"obner basis of $I$ and are real algebraic in nature, as we never compute any complex zero.

\bigskip\noindent
{\bf \em Related literature.}
There is a large literature on the problem of finding the radical ideal $\sqrt I$ of an ideal $I$; see, e.g., 
\cite{BW96}, \cite{CCT97}, \cite{EHV92}, \cite{GTZ88},  \cite{krick}, \cite{LL91}.
For the general (positive-dimensional) case, Krick and Logar \cite{krick} propose an efficient algorithm based on splitting and reduction to the zero-dimensional case, which is implemented e.g. in the computer algebra package Singular \cite{GPS05}.
In the zero-dimensional case the problem is considered to be well-solved, e.g., via the following method of Seidenberg~\cite{Se74}:
$\sqrt I= \langle I\cup\{q_1,\ldots,q_n\}\rangle$, the ideal generated by $I$ and the $q_i$'s,
where $q_i$ is the square-free part of the monic
generator $p_i$ of $I\cap\oK[x_i]$. Finding $p_i$ is easy once a linear basis $\BB$ of 
$\PK/I$ is known. Namely, find the smallest integer $k_i$ for which 
$\{1,x_i,x_i^2,\ldots,x_i^{k_i}\}$ is linearly dependent in $\PK/I$; then this smallest linear dependence gives the polynomial $p_i$.
Next, the polynomial $q_i$ can be found taking derivatives and gcd-computations as
$q_i=\frac{p_i}{\text{gcd}(p_i,p_i')}$.
 So finding $I(\VC(I))$ is easy 
if we have a basis of $\PK/I$. A classical method for finding such a
 basis $\BB$ is to compute a Gr\"obner basis of $I$ 
and the corresponding set $\BB$ of standard monomials.
The results in the present paper show that one may alternatively find such a
 basis $\BB$ from a suitable moment matrix.

On the other hand, the problem of computing the real radical ideal is considered to be much more difficult. 
For instance, in their paper  which is one of the first classical references on this problem,
Becker and Neuhaus \cite[p. 7]{BN93} write that {\em the computation of $\tau$-real parts} (thus, the real radical ideal) {\em  is much more difficult} (than that of the ordinary radical).
They give an algorithm  for $\sqrt[\oR]{I} $ based on finding the minimal real prime ideals $P_i$ such that
$\sqrt[\oR]{I}=\cap_iP_i$. Among other advanced algebraic manipulations,  their algorithm makes intensively use of 
(ordinary) radical computations. 
For other works along similar lines see, e.g., 
\cite{BS97}, \cite{CT98}.

Finally, let us mention that excellent algorithms and software packages
exist for computing the complex variety
$\VC(I)$ of a zero-dimensional ideal $I$, e.g., by 
Verschelde \cite{verschelde}, 
 Rouillier \cite{rouillier};
see also
related work by   Mourrain
et al. \cite{mourrain} and e.g. the monograph \cite{DiEm05}.
For instance, Verschelde \cite{verschelde} proposes symbolic-numeric algorithms via homotopy continuation methods (cf. also \cite{SW05}) whereas Rouillier \cite{rouillier} solves a zero-dimensional system of polynomials by giving a 
rational univariate representation (RUR) for its solutions, of the form $f(t)=0$, $x_1=\frac{g_1(t)}{g(t)}$, $\ldots$, $x_n=\frac{g_n(t)}{g(t)}$, where $f,g,g_1,\ldots,g_n\in \oK[t]$ are univariate polynomials. The computation of the RUR relies in an essential way on the multiplication matrices in the quotient algebra $\PK/I$ which thus requires the knowledge of a corresponding linear basis. 

\bigskip\noindent
{\bf \em Our contribution.} 
Given an ideal $I\subseteq \oK[x]$ ($\oK=\oR,\oC$) defined by a set of generators and satisfying $|V_\oK(I)|<\infty$, we provide a method for computing $V_\oK(I)$
as well as a border basis and a Gr\"obner basis of 
the ideal $I(V_\oK(I))$. Our approach is based on a semidefinite programming characterization of $I(V_\oK(I))$ with the following distinguished feature.
Remarkably, {all} information needed to compute the above objects is contained in the so-called 
{moment matrix} (whose entries depend on the polynomials generating the ideal $I$) 
and the geometry behind it when this matrix is required to be positive semidefinite with maximum rank. The latter property is achieved by standard semidefinite programming algorithms. 
For the task of computing the real roots and the real radical ideal $I(V_\oR(I))$, the method is real algebraic in nature, as
we do not compute (implicitly or explicitly) any complex element of $\VC(I)$.

Lasserre \cite{Las01} already recognized that moment matrices can be used 
for approximating the minimum of a polynomial over a basic closed semi-algebraic set and sometimes extracting global minimizers (cf. \cite{HeLa05}).
The present  paper builds on this approach and shows how it can be applied to finding the real radical of a zero-dimensional ideal.
 Moreover there are links between moment matrices and the 
 Hermite quadratic forms used in \cite{PRS93} for computing the number of real roots, that were pointed out in \cite{Lau1}.

Our approach with its specificity is best illustrated  on the task
of computing the real radical ideal $I(V_\oR(I))$ that we now briefly describe.

\medskip
Given a sequence $y=(y_\alpha)_{\alpha\in \oN^n}\in\oR^{\oN^n}$, consider the 
{\em moment matrix}
$$M^\oR(y):=(y_{\alpha+\beta})_{\alpha,\beta\in\oN^n}$$
(later we will also introduce complex moment matrices $M^\oC(y)$, $\fCM(y)$).
One may think that $y$ and $M^\oR(y)$ are   indexed by the set $\oT_n:=\{x^\alpha\mid
\alpha\in\oN^n\}$ of monomials.
Given a polynomial $h\in \PR$, set $\ve{h}:=(h_\alpha)_{\alpha\in\oN^n}$ and 
define the new sequence $hy:=M^\oR(y)\ve{h} \in \oR^{\oN^n}$. By abuse of language let us say that $h$ lies in the kernel of $M^\oR(y)$ when $\ve{h}$ does, which enables us to view $\Ker M^\oR(y)$ as a subset of $\PR$.
The following property of moment matrices plays a central role
in our approach; it is based on ideas from \cite{CF96},\cite{CF00},\cite{Lau05} and will be proved at the end of  Section \ref{secprelmom}.
Let $I=\langle h_1,\ldots,h_m\rangle$ be an  ideal generated by $h_1,\ldots,h_m\in \oR[x]$.

\begin{proposition}\label{prop0}
Assume that $V_\oR(I)$ is finite. If 
\begin{equation}\label{eq0}
M^\oR(y)\succeq 0,\ M^\oR(h_jy)=0 \ (j=1,\ldots,m)
\end{equation}
then  the kernel of $M^\oR(y)$ is a real radical ideal,
$\rank M^\oR(y) \le |V_\oR(I)|$
and 
$I(V_\oR(I))\subseteq \Ker M^\oR(y)$, with equality if and only if 
$M^\oR(y)$ has maximum rank, equal to $|V_\oR(I)|$.
\end{proposition}

\noindent
(In (\ref{prop0}) the notation "$\succeq0$" stands for positive semidefinite.)
This semidefinite characterization leads directly to an algorithm for computing
$I(V_\oR(I))$, by considering {\em truncated} moment matrices in place of the full (infinite) moment
matrix $M^\oR(y)$. Namely, given an integer $t$, let $M^\oR_t(y)$ denote the 
principal submatrix of $M^\oR(y)$ whose rows and columns are indexed by the set
$\oT_{n,t}:=\{x^\alpha\mid \alpha \in\oN^n \ \text{ with }   |\alpha|:=\sum_i\alpha_i\le t\}$
 and set 
 \begin{equation}\label{reldj}
d_j:=\lceil \deg(h_j)/2\rceil, \ d:=\max_{j=1,\ldots,m}d_j.
\end{equation}
Fix $t\ge d$ and assume $M^\oR_t(y)$ is a maximum rank matrix satisfying
\begin{equation}\label{mat0}
M^\oR_t(y)\succeq 0, \ M^\oR_{t-d_j}(h_jy)=0 \ (j=1,\ldots,m).
\end{equation}
We will show that if, moreover,  
\begin{equation}\label{rankd}
\rank M^\oR_s(y)=\rank M^\oR_{s-d}(y)
\end{equation}
for some $d\le s\le t$, then
$I(V_\oR(I))$ coincides with the ideal generated by $\Ker M^\oR_s(y)$.
The same conclusion holds if
\begin{equation}\label{rank1}
\rank M^\oR_s(y)=\rank M^\oR_{s-1}(y)\end{equation}
 for some  $2d\le s\le t$.
Moreover, from the semidefinite characterizations
(\ref{mat0})-(\ref{rank1}), the following algebraic objects
can be obtained directly from the matrix $M^\oR_t(y)$:

\begin{description}
\item[(i)]Let  $\BB\subseteq \oT_{n,s}$ be a set indexing a maximum nonsingular principal submatrix of $M^\oR_s(y)$. 
Then  $\BB$ 
is a linear basis of the quotient vector space $\PR/I(V_\oR(I))$ (see Section \ref{secflat}).
\item[(ii)]  We can compute directly from $M^\oR_t(y)$ 
the matrix of any  multiplication operator in $\PR/I(V_\oR(I))$ with respect to the basis $\BB$,  
and thus compute $V_\oR(I)$ (using the eigenvalue method, see Section \ref{seceig}).
\item[(iii)] 
 When the set $\BB$ (as in (i)) is an order ideal (i.e., is stable under division),
 the matrices of the multiplication operators by $x_1,\ldots,x_n$ give directly a border basis of the ideal $I(V_\oR(I))$ (see Section \ref{secborder}).
\item[(iv)] Given a graded lexicographic monomial ordering, we can find a set $\BB$ (as in (i))
 which is precisely the set of standard monomials; the associated reduced Gr\"obner basis of $I(V_\oR(I))$ can then be recovered, since it is contained in the border basis. In fact our method also applies to an arbitrary monomial ordering (see Section \ref{secgrobner}).
\item[(v)] Finally the method can also detect whether the real variety $V_\oR(I)$ is empty. Indeed, $V_\oR(I)=\emptyset$ if and only if,
for some integer $t$, the system (\ref{mat0}) admits no solution $y$ with $y_0\ne 0$. (See Remark \ref{remdetect}.)
\end{description}
\medskip
\noindent
{\bf \em Further discussion.} 
An independence oracle in $\PR/I(V_\oR(I))$ is needed for our algorithm in (iv) above. The following property is a 
crucial ingredient. Assume 
that one of the conditions (\ref{rankd}) or (\ref{rank1}) holds and consider a set 
$T\subseteq \oT_{n,s}$. Then, $T$ is linearly independent 
 in $\PR/I(V_\oR(I))$ if and only if  $T$ indexes a linearly independent set 
of columns 
 of $M^\oR_t(y)$.
In view of (iv),  a Gr\"obner basis can easily be derived {\em afterwards} in contrast with
 classical methods which compute the set of standard monomials from the Gr\"obner basis.

Realizing the above tasks relies only on numerical linear algebraic operations on $M^\oR_t(y)$ like 
evaluating the rank of certain principal submatrices.  
Finding a matrix satisfying (\ref{mat0}) is an instance of  semidefinite programming. Moreover, 
it is a property of most interior-point algorithms for semidefinite programming that they do find such a matrix having maximum rank (see Section~\ref{rem::SDP} for details).

 The method is iterative. Namely, if the maximum rank matrix satisfying
(\ref{mat0}) does not satisfy (\ref{rankd}) or (\ref{rank1}), then 
we iterate with $t+1$ in place of $t$.  
The method eventually terminates since we will  show that (\ref{rankd}) holds for $t$ large enough.

The following two small examples illustrate 
how positive semidefiniteness of the matrix $M^\oR_t(y)$
allows the elimination of  all complex (nonreal) roots, whose number can be much larger than the number of real roots or even infinite. 

\begin{example}
\label{exm::jean}
Let $I\subseteq \oR[x]$ be generated by $h_i=x_i(x_i^2+1)$ ($i=1,..,n$).
Then, $V_\oR(I)=\{0\}$,  $|\VC(I)|=3^n$, $d_i=2$ for all $i$.
Assume  $y$ satisfies (\ref{mat0}) for order $t=3$. 
Then $M^\oR_1(h_iy)=0$ implies $y_{4e_i}=-y_{2e_i}$
and $M^\oR_3(y)\succeq 0$ implies  $y_{2e_i},y_{4e_i}\ge 0$ which in turn
implies $y_{\alpha}=0$ for all $\alpha\ne 0$ with $|\alpha| \le 5$.
(Throughout, $e_1,\ldots,e_n$ denote the standard unit vectors in $\oR^n$.)
Hence $\rank M^\oR_2(y)= \rank M^\oR_0(y)=1$; that is,
  (\ref{rankd}) holds for $s=2$. 
In fact, $\Ker M^\oR_1(y)$ is spanned by $x_1,\ldots,x_n$, the generators 
of $I(V_\oR(I))$. One may argue that the ideal $I$ is already described by a Gr\"obner basis. But the same conclusion 
also holds under the change of variables $x=Ay$ with $A$ being a nonsingular matrix, in which case other methods would require a Gr\"obner basis computation.
\end{example}

\begin{example}\label{ex2}
Let $I\subseteq \oR[x_1,x_2]$ be generated by $h=x_1^2+x_2^2$.
Then $V_\oR(I)=\{0\}$ and $\VC(I)=\{(x_1,x_2)\mid x_1=\pm ix_2\}$ is infinite.
Then  $M^\oR_0(hy)=0$ gives $y_{2e_1}+y_{2e_2}=0$ which, together with $M^\oR_1(y)\succeq 0$,  implies 
$y_\alpha=0$ for $\alpha\ne 0$. Hence the maximum rank of $M^\oR_1(y)$ is equal to 1 and 
$\Ker M^\oR_1(y)$ is spanned by $x_1,x_2$, the generators of $I(V_\oR(I))$.
\end{example}

The method sometimes (partially) applies even if none of the rank conditions 
(\ref{rankd}), (\ref{rank1}) holds, namely
when $M^\oR_t(y)$ contains sufficient information for the construction  of
the (formal) multiplication matrices.
More precisely, let  $M^\oR_t(y)$ be a maximum rank matrix satisfying (\ref{mat0}),
let $\BB\subseteq \oT_{n,t}$ index a maximum nonsingular 
principal submatrix of $M^\oR_t(y)$,  set 
$\partial \BB:=(\cup_{i=1}^n x_i\BB ) \setminus \BB$,
and assume that the two principal submatrices
of $M^\oR_t(y)$ indexed by $\BB$ and by $\BB \cup \partial \BB$ have the same rank.
Then, by the results of  Kehrein, Kreuzer and Robbiano~\cite[Ch. 4]{DiEm05}, we can construct the formal multiplication matrices and, if they commute pairwise, a set $W\supseteq V_\oR(I)$ can be computed. By checking whether the points of $W$ satisfy all the equations $h_j=0$, we can 
eliminate the points in $W\setminus V_\oR(I)$. It turns out that,
for most examples we have tested,
$W=V_\oR(I)$ holds and we are again able to find $I(V_\oR(I)$ together with a border basis generating this ideal.

Finally, the method also applies to the task of finding the radical ideal
$I(\VC(I))$ of a zero-dimensional ideal $I\subseteq \PC$. For this, instead of using the matrix $M^\oR_t(y)$ where $y$ is a real sequence indexed by $\oT_{n,t}$, we have to use a matrix $\fCM_t(y)$ where the argument 
is a complex sequence indexed by $\oT_{2n,t}$ (see Section \ref{prelmommat} for details). Similar results hold as in the real case. Namely, under certain rank conditions, the ideal $I(\VC(I))$ can be obtained as the ideal generated by the kernel of a maximum rank complex moment matrix (see Section \ref{secfind3} for details). However, a drawback in the complex case is that one must in general handle matrices of larger order which leads to larger semidefinite programs,
thus more difficult to solve. 
However, so far we do not claim that our method can compete with existing methods
for finding the complex variety $\VC(I)$ as e.g.  \cite{rouillier}, or  
\cite{verschelde}, especially in view of the present status of SDP solvers
(that we use as a black box), still in their infancy.

\bigskip\noindent
{\bf \em Contents of the paper.}
Section \ref{secprelideal} contains preliminaries about ideals of polynomials, in particular, about
the quotient ring $\PK/I$,  multiplication matrices, Gr\"obner bases and border bases.  
We also indicate in Section \ref{secsieve} an algorithm for finding the set of 
standard monomials from an independence oracle in $\PK/I$.
Section \ref{secprelmom} contains preliminaries about moment matrices,
in particular, results relating (real) radical ideals and kernels of positive semidefinite moment matrices.
In Section \ref{secfind}, we prove the main results about the semidefinite characterization
of the variety $V_\oK(I)$ and the associated radical ideal
$I(V_\oK(I))$. Section \ref{subsec::algorithm} gives the details and implementation of an algorithm based
on the semidefinite characterization, and Section \ref{sec::numexamples}
contains several examples illustrating its behaviour.

\section{Preliminaries on Polynomial Ideals}\label{secprelideal}
\subsection{Polynomial ideals and varieties}\label{secprelideal1}
Throughout, $\oK=\oR$ or $\oC$, and $\PK:=\oK[x_1,\ldots,x_n]$ 
denotes the ring of multivariate polynomials in $n$ variables over the field $\oK$.
For an integer $t\ge 0$, $\oK[x]_t$ denotes the set of polynomials of degree at most $t$.
For a scalar $a\in \oC$, $\bar a$ denote its complex conjugate and,
 for a vector $u\in \oC^n$ (resp., a matrix $A$), $u^*$ (resp., $A^*$)
 denotes its conjugate transpose. 
Following e.g. \cite{CLO97}, $x^\alpha$ denotes the monomial $x_1^{\alpha_1}\cdots x_n^{\alpha_n}$ (for $\alpha\in\oN^n$) and $cx^\alpha$ is called a term (for $c\in\oK$).  
Let  $\oT_n:=\{x^\alpha\mid \alpha\in \oN^n\}$  denote the set of monomials and set
$\oN^n_t:=\{\alpha\in\oN^n\mid |\alpha|:=\sum_{i=1}^n\alpha_i\le t\}$,
$\oT_{n,t}:=\{x^\alpha\mid \alpha\in \oN^n_t\}$ for
$t\in \oN$. 
Following \cite[Ch.~4]{DiEm05}, a set $\BB\subseteq \oT_n$ is called an {\em order ideal}
if $\BB$ is stable under division, i.e., for all
$a,b\in\oT_n$, $b\in\BB$, $a|b$ implies $a\in \BB$.
Given an ideal\footnote{An ideal $I\subseteq\PK$ is an additive subgroup of 
$\PK$ such that $fg\in I$ whenever $f\in I$ and $g\in\PK$.}  $I\subseteq \PK$, let
$$V_\oC(I):=\{x\in\oC^n\mid f(x)=0\ \forall f\in I\}, \
V_\oR(I):=\VC(I)\cap \oR^n$$
denote its complex and real varieties, respectively.
For a set $V\subseteq \oK^n$,  define the ideal 
$$I(V):=\{f\in \PK\mid f(v)=0\ \forall v\in V\}.$$
Given an ideal $I\subseteq \PK$, one can define the ideals
$I(\VC(I))$ and
$$\sqrt I:=\{f\in \PK\mid f^m\in I \ \text{ for some } m\in \oN\setminus \{0\}\}$$
and, when $I\subseteq \PR$, one can define the  ideals 
$I(V_\oR(I))$ and
$$\sqrt[\oR]I:=\{p\in \oR[x]\mid p^{2m} +\sum_j q_j^2\in I \ \text{ for some }q_j\in \oR[x], m\in \oN\setminus \{0\}\}.$$ 
Obviously, 
$$I\subseteq \sqrt I\subseteq I(\VC(I)),\ \   I\subseteq \sqrt[\oR]I\subseteq I(V_\oR(I)).$$ 
The ideal $I$ is said to be {\em radical} (resp., {\em real radical}) if $I=\sqrt I$ (resp., $I= \sqrt[\oR]I$).
Obviously, $I\subseteq I(\VC(I))\subseteq I(V_\oR(I))$. Hence, if $I\subseteq \PR$ is real radical, then 
$I$ is radical and moreover,
$\VC(I)=V_\oR(I)\subseteq \oR^n$ if $|V_\oR(I)|<\infty$.
The following lemma gives a useful characterization for  (real) radical ideals.

\begin{lemma} \label{lemradical}
An ideal $I\subseteq \PK$ is radical if and only if 
\begin{equation}\label{radical}
\forall  p\in\PK, \ \ \  p^2\in I\Longrightarrow p\in I 
\end{equation}
and $I\subseteq \PR$ is real radical if and only if 
\begin{equation} \label{realradical}
\forall p_1,\ldots, p_k\in\PR, \ \ \ 
p_1^2+\ldots +p_k^2 \in I \Longrightarrow p_1,\ldots,p_k \in I .
\end{equation}
\end{lemma}

\begin{proof}
If $I$ is radical, then (\ref{radical}) obviously holds.
Conversely assume that
(\ref{radical}) holds; we show that 
$p^m\in I$ $\Longrightarrow  p\in I$ by induction on 
$m\ge 2$. If $p^m\in I$ then $p^{m+1}\in I$ and thus, by 
 (\ref{radical}), 
$p^{\lceil m/2\rceil}\in I$ which, by the induction assumption, implies $p\in I$.
The proof for the real radical case is along the same lines and thus omitted.
$\hfill\Box$
\end{proof}

\begin{theorem}
\begin{description}
\item[(i)] {\bf Hilbert's  Nullstellensatz} (see, e.g., \cite[\S 4.1]{CLO97}) $\sqrt I=I(\VC(I))$.
\item[(ii)] {\bf Real Nullstellensatz} (see, e.g., \cite[\S 4.1]{BCR})
$ \sqrt[\oR]I=I(V_\oR(I))$ for  an  ideal $I\subseteq \PR$.
\end{description}\end{theorem}

For a polynomial $p\in \oR[x]$, $x\mapsto p(x)=\sum_\alpha p_\alpha x^\alpha$, let
$\ve{p}:=(p_\alpha)_{\alpha\in\oN^n}$ denote the vector of its
coefficients. We also let $\ve{p}$ denote the vector
$(p_\alpha)_{\alpha\in\oN^n_t}$ for any $t\ge \deg(p)$, as $p_\alpha=0$
whenever $|\alpha|>\deg(p)$.

Finally, given $A\subseteq \oK[x]$, let $\langle A\rangle:=\{\sum_{i=1}^ku_ip_i\mid
p_i\in A,  u_i\in \oK[x]\}$ denote the ideal generated by $A$.

\subsection{The algebra $\PK/I$ and multiplication matrices}\label{seceig}

Consider the quotient space $\oK[x]/I$, whose elements are the cosets
$[f]:=f+I=\{f+q \mid q \,\in\,I\}$ for $f\in \PK$. 
$\oK[x]/I$  is a $\oK$-vector space with addition
$ [f]+[g]:=[f+g]$ and scalar multiplication $\lambda [f]:=[\lambda f]$,
and an algebra with multiplication
$[f][g]:=[fg]$, 
 for $\lambda\in \oK$, $f,g\in \PK$.
In particular, for $h\in \PK$, the {\em multiplication operator}
$$\begin{array}{lccc}
m_h: & \PK/I & \longrightarrow & \PK/I\\
     &  [f] & \longmapsto & [hf]
\end{array}$$
is well defined. 
The following  well known result 
relates  the cardinality of $\VC(I)$ and
the dimension of the vector space $\PK/I$. See e.g.  \cite{CLO97}, \cite{St04} for a detailed treatment of the quotient algebra $\PK/I$.

\begin{theorem}\label{theodim}
For an ideal $I$ in $\PK$,
$|\VC(I)|<\infty \Longleftrightarrow
\dim \PK/I<\infty.$
Moreover, $|\VC(I)|\le \dim\ \PK/I$, with equality if and only if  $I$ is radical.
\end{theorem}

Assume $|\VC(I)|<\infty$ and set $N:=\dim \PK/I\ge |\VC(I)|$. Consider a set 
 $\BB:=\{b_1,\ldots,b_N\}\subseteq \PK$ for
which the cosets $[b_1],\ldots,[b_N]$ are pairwise distinct and 
$\{[b_1],\ldots,[b_N]\}$ is a basis
of $\PK/I$; by abuse of language we also say that $\BB$ itself is a basis of
$\PK/I$.
Then  every $f\in \PK$ can be written in a unique way as
$f=\sum_{i=1}^N c_i b_i +p,$  where  $c_i\in \oK,$ $ p\in I;$
the polynomial $\res_\BB(f):=\sum_{i=1}^N c_i b_i$ is called the {\em residue of $f$ modulo $I$ w.r.t. the basis $\BB$}.
In other words, the vector space $\Span_\oK(\BB):=\{\sum_{i=1}^N c_ib_i\mid c_i\in \oK\}$ is isomorphic to
$\PK/I$.

Given a basis $\BB$ of $\PK/I$ and $h\in\PK$, let $\mathcal{X}_h$ denote the matrix of the multiplication 
operator $m_h$ with respect to $\BB$. That is, writing $\res_\BB(hb_j)=
 \sum_{i=1}^N a_{ij}b_i$, the $j$th column of $\mathcal{X}_h$ is the vector
$(a_{ij})_{i=1}^N$.
The following well known result relates the points of the variety $\VC(I)$ to the eigenvalues and eigenvectors of $\mathcal{X}_h$. See, e.g.,  \cite[Ch. 2,3]{DiEm05} for a detailed treatment.

\begin{theorem} \label{theoeig}
Let $h\in \PK$ and, for $v\in \VC(I)$, set 
$\zeta_{\BB,v}:=(b_i(v))_{i=1}^N$.
The set $\{h(v)\mid v\in \VC(I)\}$ is the set of eigenvalues of $\mathcal{X}_h$ and
 $\mathcal{X}_h^T\zeta_{\BB,v}=h(v)\zeta_{\BB,v}$ for all
$v\in \VC(I)$. 
\end{theorem}

\noindent
When the matrix $\mathcal{X}_h$ is non-derogatory (i.e., all its eigenspaces are
 1-dimensional),  one can recover 
the points $v\in \VC(I)$ from the eigenvectors of $\mathcal{X}_h^T$. If
$I$ is radical, then $N=|\VC(I)|$ and thus 
$\mathcal{X}_h$ is non-derogatory whenever the values $h(v)$ ($v\in \VC(I)$) are pairwise distinct. 
This is achieved with high probability if one
chooses $h=\sum_{i=1}^n a_ix_i$ for random 
scalars $a_i$.

\subsection{Gr\"obner bases and standard monomials}
A classical basis  of $\PK/I$ is 
the set of standard monomials 
with respect to some monomial ordering 
 `$\succ $' of $\oT_n$. 
Let us recall some definitions. (See e.g. \cite{CLO97} for details.)
Fix  a monomial ordering $\succ$ on $\oT_n$. Write also
$ax^\alpha \succ bx^\beta$ if $x^\alpha \succ x^\beta$ and $a,b\in \oK\setminus  \{0\}$.
For a nonzero polynomial $f=\sum_\alpha f_\alpha x^\alpha$,
its {\em leading term} $\LT(f)$ is the maximum $f_\alpha x^\alpha$ with respect to $\succ$ for which $f_\alpha\ne 0$. 
The leading term ideal of $I$ is
$\LT(I):=\langle \LT(f)\mid f\in I\rangle$ and 
the set
$$\BB_\succ:=\oT_n\setminus \LT(I)=\{x^\alpha\mid \LT(f) \text{ does not divide } x^\alpha \ \  \forall f\in I\}$$
is the set of {\em standard monomials}. Obviously  $\BB_\succ$ is an order ideal.
 A finite set $G\subseteq I$ is a {\em Gr\"obner basis} of $I$
 if $\LT(I)=\langle \LT(g)\mid g\in G\rangle$; thus $x^\alpha\in\BB_\succ$  
if and only if $x^\alpha$  is not divisible by the leading term of any
 polynomial in $G$. 
A Gr\"obner basis always exists and it can be constructed, e.g., with the algorithm of Buchberger. 
Call $G$ {\em reduced} 
if, for all $g\in G$, the leading coefficient of $\LT(g)$ is 1 and no 
term of $g$ lies in
 $\langle         LT(g') \mid g'\in G\setminus\{g\}\rangle$.
Given nonzero polynomials $f,h_1,\ldots,h_m$,
 the division algorithm applied to dividing 
  $f$ by $h_1,\ldots,h_m$ produces 
 polynomials $u_1,\ldots,u_m,r$ satisfying
  $f=\sum_{j=1}^m u_jh_j +r,$  no term of $r$ is divisible by
  $\LT(h_j)$ ($j=1,\ldots,m$)  and $ \LT(f) \succeq  \LT(u_jh_j)$.
Note that $\deg(u_ih_i)\le \deg(f)$ when the monomial ordering is a
 graded lexicographic  ordering.
 When  $\{h_1,\ldots,h_m\}$ is a Gr\"obner basis 
of the ideal $I:=\langle h_1,\ldots,h_m\rangle$, the  remainder $r$ is uniquely determined and 
belongs to $\Span_\oK(\BB_\succ)$; moreover, $f\in I\Longleftrightarrow r=0$.
Therefore, the set $\BB_\succ$ is a basis of $\PK/I$.

\medskip
For an arbitrary  basis $\BB$ of $\PK/I$, 
set $d_\BB:= \max_{b\in \BB}  \deg(b)$. The next result shows that $d_\BB$ is minimum when $\BB$
is the set of standard monomials for some graded lexicographic order.

\begin{lemma}\label{lemgrad}
Let $I$ be a zero-dimensional ideal in $\PK$. 
Let $\{g_1,\ldots,g_k\}$ be the Gr\"obner basis of $I$ with respect to 
a graded lexicographic monomial ordering and let $\BB_\succ$ be 
the corresponding set of standard monomials.
For any basis $\BB$ of $\oK[x]/I$, we have $d_{\BB_\succ}\le d_{\BB}.$
\end{lemma}

\begin{proof}
Set $\BB=\{b_1,\ldots,b_N\}$. Write
$b_i=\sum_{x^\alpha\in \BB_\succ} c_{i,\alpha}x^\alpha +
\sum_{h=1}^k u_hg_h$
where $ c_{i,\alpha}\in \oK$ for $i=1,\ldots,N$ and $u_h\in\PK$. Then $\deg(u_hg_h)\le \deg(b_i)$ (by the properties of the division algorithm as we use a graded monomial ordering). Thus, 
$\deg(\sum_{x^\alpha\in \BB_\succ} c_{i,\alpha}x^\alpha)\le \deg(b_i)$.
Let $x^{\alpha_0}\in\BB_\succ$ with $\deg(x^{\alpha_0})=d_{\BB_\succ}$. 
As $\BB$, $\BB_\succ$ are two bases of $\oK[x]/I$, the 
matrix $(c_{i,\alpha})_{\genfrac{}{}{0pt}{1}{i=1,\ldots,N}{\alpha\in\BB_\succ}}$
is nonsingular and thus its $\alpha_0$th column is nonzero. Hence $c_{i,\alpha_0}\ne 0$ for some $i$.
Hence $d_{\BB_\succ}=\deg(\sum_{x^\alpha\in \BB_\succ} c_{i,\alpha}x^\alpha)
\le \deg(b_i)\le d_\BB$.
$\hfill\Box$
\end{proof}

\subsection{Finding the set of standard monomials from an independence oracle}\label{secsieve}
When $I$ is a zero-dimensional ideal and $\succ$ is a monomial ordering on $\oT_n$, we describe a method for finding the set $\BB_\succ$ of standard monomials, assuming we have an oracle 
for checking  linear independence in $\PK/I$. This \emph{`greedy sieve' algorithm}, described below in Algorithm~\ref{alg::greedysieve}, does not require  knowledge of a Gr\"obner basis of $I$. 

\begin{algorithm}
\caption{\emph{The `greedy sieve' algorithm:}}
\begin{algorithmic}[1]
\Require A zero-dimensional ideal $I\in \PK$, a monomial ordering $\succ$ on $\oT_n$, and an integer $s\ge 1$.
\Ensure A set $\BB\subseteq \oT_{n,s}$ linearly independent in 
$\PK/I$ and satisfying $\BB\supseteq \BB_\succ\cap \oT_{n,s}$
\State Order  the monomials in $\oT_{n,s}$ with respect to $\succ$.
\State Initialize $\BB:=\emptyset$, $L:=(t_1,t_2,\ldots)$, the ordered set $\oT_{n,s}$.
\While{$\BB\subset L$}
\State Set $t$ as the first element of $L\setminus \BB$
\If{$\BB\cup\{t\}$ is linearly independent in $\PK/I$}
\State Reset $\BB:=\BB\cup\{t\}$
\Else
\State Reset $L:=L\setminus t\oT_n$ (i.e., remove from $L$ all multiples of $t$).
\EndIf
\EndWhile
\State \Return $\BB=L$
\end{algorithmic}
\label{alg::greedysieve}
\end{algorithm}

The next lemma shows correctness of Algorithm 1 and how to use it for finding the set $\BB_\succ$ of standard monomials.

\begin{lemma} \label{lemalg}
Let $I$ be a zero-dimensional ideal,  $\succ$ a monomial ordering on $\oT_n$, and 
$\BB_\succ=\oT_n\setminus \LT(I)$, the associated set of standard monomials.
For an integer $s\ge 1$, let $\BB_s$ be the set returned by 
the greedy sieve algorithm applied to $(I,\succ,s)$. 
\begin{description}
\item[(i)]
$\BB_s$ is linearly independent in $\PK/I$ and satisfies $\BB_\succ\cap\oT_{n,s}\subseteq \BB_s$;
in particular, $\BB_s=\BB_\succ$  if $\BB_\succ\subseteq \oT_{n,s}$.
\item[(ii)] If $\BB_s=\BB_{s+1}$ then $\BB_s=\BB_\succ$.
\item[(iii)]
If $\succ$ is a graded monomial ordering, then $\BB_s\subseteq \BB_\succ$; therefore,
$\BB_s=\BB_\succ$ if $|\BB_s|=\dim \PK/I$.
\end{description}
\end{lemma}

\begin{proof}
(i) Obviously, throughout the algorithm, $\BB$ is linearly independent in $\PK/I$ and
$\BB\subseteq L$.
Assume  $t_k\in (\BB_\succ\cap \oT_{n,s})\setminus \BB_s$.
Consider the step when the algorithm examines $t_k$ and let $\BB$ be the current set maintained by the algorithm.
Then, $\BB\subseteq \{t_1,\ldots,t_{k-1}\}$,
$t_k\in L\setminus \BB$ and $\BB\cup\{t_k\}$ is linearly dependent in $\PK/I$.
Hence there exists a polynomial $f\in I$ with $\LT(f)=t_k$,
contradicting the assumption that $t_k\in \BB_\succ$. This shows $\BB_\succ\cap\oT_{n,s}\subseteq \BB_s$.
Moreover, if $\BB_\succ\subseteq \oT_{n,s}$, then $\BB_\succ \subseteq \BB_s$;
equality holds since $|\BB_s|\le \dim\PK/I$ as $\BB_s$ is linearly independent in
$\PK/I$, while $|\BB_\succ|=\dim\PK/I$. \\
(ii) 
Assume $\BB_s=\BB_{s+1}$. Then, in view of (i), $\BB_\succ \cap (\oT_{n,s+1}\setminus \oT_{n,s})=\emptyset$. This 
implies $\BB_\succ \subseteq \oT_{n,s}$. 
Indeed assume $t\in \BB_\succ$ has degree at least $s+1$; then any divisor $t'$
of $t$ with degree $s+1$ lies in $\BB_\succ$ (since $\BB_\succ$ is an order ideal) and thus $t'\in \BB_\succ \cap (\oT_{n,s+1}\setminus \oT_{n,s})=\emptyset$, a contradiction.
Therefore, by (i), $\BB_s=\BB_\succ$.\\
(iii)
Assume $\succ$ is a graded monomial ordering and, say, $\BB_\succ\subseteq \oT_{n,d}$ for some integer $d$.
If $d\le s$ then $\BB_s=\BB_\succ$ by (i).
If $d\ge s+1$, then $\BB_s\subseteq \BB_d$  
(since all elements of $\oT_{n,d}\setminus \oT_{n,s}$ 
come after the elements of $\oT_{n,s}$ in the ordering $\succ$) and $\BB_d=\BB_\succ$ (by (i)), implying
$\BB_s\subseteq \BB_\succ$.
$\hfill\Box$
\end{proof}

\begin{remark}\label{remsieve}
Observe that, when $\succ$ is not a graded monomial degree ordering, one cannot claim the inclusion $\BB_s\subseteq \BB_\succ$. For instance, consider the ideal
$I=\langle x^3-1, -y+x^2+x+1\rangle$ in $\oR[x,y]$ and choose as monomial ordering $\succ$ the lexicographic order with $y>x$. 
Then, $\BB_\succ=\{1,x,x^2\}$ is found when applying Algorithm 1 to 
$(I,\succ,s=2)$; observe that $\BB_2=\BB_3$. However,
the algorithm applied to $(I,\succ,s=1)$ returns the set
$\BB_1=\{1,x,y\}$;  thus $\BB_1\not\subseteq \BB_\succ$, while
$|\BB_1|=3=\dim\PR/I$.

Alternatively, one could  initialize the set $L$ in Algorithm 1 to be the full ordered set $\oT_n$.
Then the algorithm still terminates in finitely many steps (because $\dim\PK/I<\infty$) and the set $\BB$ 
returned by the algorithm is equal to $\BB_\succ$ (using the same argument as in Lemma \ref{lemalg} (i), one can show that $\BB_\succ \subseteq \BB$, implying $\BB_\succ=\BB$).

A crucial tool for applying Algorithm 1 is having an oracle for testing linear independence in 
$\PK/I$. 
In our setting the oracle will work as follows:
Given a subset $\BB\subseteq \oT_{n,s}$, $\BB$ is linearly independent in 
$\PK/I$ if and only if $\BB$ indexes a linearly independent set of columns of a suitable moment matrix $M_s(y)$. 
This motivates why in our presentation of Algorithm 1 we explore the set $\oT_{n,s}$ of monomials of degree at most $s$.
\end{remark}

\subsection{Border bases and formal multiplication matrices}\label{secborder}

We recall results about border bases following the exposition from \cite[Ch.~4]{DiEm05}. See also \cite{St04} 
for details about border bases. 
Given an order ideal $\BB\subseteq \oT_n$, the {\em border} of $\BB$ is the set 
\begin{equation}\label{relborder}
\pBB:=\{x_ix^\beta\mid x^\beta\in\BB,\ i=1,\ldots,n\}\setminus \BB.
\end{equation}
Assume $\BB\ne \emptyset$, set $N:=|\BB|$, $H:=|\pBB|$ and write
 $\BB=\{b_1,\ldots,b_N\}$ and
$\pBB=\{c_1,\ldots,c_H\}$.
A set of polynomials $G=\{g_1,\ldots,g_H\}$ is called a {\em $\BB$-border prebasis}  
if each $g_j$ is of the form 
\begin{equation}\label{relpreborder}
g_j=c_j -\sum_{i=1}^Na_{ij} b_i \ \text{ for some }a_{ij}\in \oK.
\end{equation}  One also says that {\em $g_j$ is  marked by the element
$c_j$} of $\pBB$.
Given a polynomial $f$, the border division algorithm \cite[Prop. 4.2.10]{DiEm05} 
produces polynomials 
$u_j,r$ such that 
$ f=\sum_{j=1}^H u_jg_j + r$, and
$r\in \Span_\oK(\BB)$.
Hence, for any ideal $I$ containing $G$,
$\BB$ spans the $\oK$-vector space $\PK/I$.
The set $G\subseteq I $ is said to be a {\em $\BB$-border basis of $I$}  if 
$\BB$ is linearly independent in $\PK/I$, i.e.,
if $\BB$ is a linear basis of $\PK/I$; in that case $G$ generates the ideal $I$.

Stetter \cite{St04} advocates using border bases instead of  Gr\"obner bases
since they do not depend on any monomial ordering. Border 
bases  represent in fact an extension of the notion of Gr\"obner bases.
Indeed, the set $\oT_n\setminus \BB$ defines a monomial ideal;  the elements of the 
minimal set of generators of this monomial ideal are called 
the {\em corners} of $\BB$, which  belong to $\pBB$.
When $\BB=\BB_\succ$ is the set of standard monomials for some monomial ordering,
there exists a unique $\BB_\succ$-border basis $G$ of $I$ and
the reduced Gr\"obner basis of $I$ is the subset of $G$ consisting 
of the polynomials in $G$ that are marked by the corners 
of $\BB_\succ$.

\medskip
When $G$ is a $\BB$-border prebasis, one can mimic the construction 
of the multiplication matrices from the previous section in the following way.
Fix $k\in\{1,\ldots,n\}$. The {\em formal multiplication matrix} 
$\XX_k$ is the $N\times N$ matrix whose $i$th column is defined 
as follows.
If $x_kb_i \in \BB$, say, $x_kb_i=b_r$, then the $i$th column of $\XX_k$ is the standard 
unit vector $e_r$ 
(with all zero entries except 1 at the $r$th position).
Otherwise, $x_kb_i\in \pBB$, say, $x_kb_i =c_j$, then the $i$th column 
of $\XX_k$ is the vector $(a_{ij})_{i=1}^N$ (compare with Eqn.(\ref{relpreborder})).
We will use the following result (see \cite[Thm. 4.3.17]{DiEm05}).

\begin{theorem}\label{theoborder} 
Let $\BB\subseteq\oT_n$ be an order ideal, let $G$ be a $\BB$-border prebasis with associated formal multiplication matrices $\XX_1,\ldots,\XX_n$,
and let $J:=\langle G\rangle$ be the ideal generated by $G$. Then,
$G$ is a border basis of $J$ 
if and only if the matrices $\XX_1,\ldots,\XX_n$ commute pairwise.
In that case, $\BB$ is a linear basis of $\PK/J$ and 
 the matrix $\XX_k$  represents the multiplication operator
$m_{x_k}$ of $\PK/J$ with respect to the basis $\BB$. 
\end{theorem}

\begin{remark}\label{remMourrain}
Following Mourrain \cite{Mou99}, call  $\BB\subseteq \oT_n$ {\em connected to 1}
if $1\in\BB$ and 
 any monomial in $\BB$ is of the form
$x_{i_1}x_{i_2}\cdots x_{i_k}$ with $x_{i_1},$ $x_{i_1}x_{i_2},$ $\ldots,$ $x_{i_1}x_{i_2}\cdots x_{i_k}~\in~\BB$.
Obviously if $\BB$ is an order ideal then $\BB$ is connected to 1.
As shown by Mourrain \cite[Th. 3.1]{Mou99},
the result of Theorem \ref{theoborder} remains valid in the more general  setting where $\BB$ is connected to 1 (instead of being an order ideal).
We restrict our attention in this paper to monomial bases of $\PK/J$ that are order ideals,
in particular, because we have an algorithm for finding such bases, as we just saw in the preceding section.
It will be interesting to investigate the use of bases satisfying Mourrain's criterion in subsequent work.
\end{remark}

\section{Preliminaries on Moment Matrices} \label{secprelmom}

\subsection{Moment matrices} \label{prelmommat}

Given  a sequence $y\in \oR^{\oN^n}$, its {\em real moment matrix}
$M^\oR(y)$ is the real symmetric matrix indexed by $\oN^n$ whose $(\alpha,\beta)$th entry
is $y_{\alpha+\beta}$, for $\alpha,\beta\in\oN^n$.
Given a sequence $y\in \oC^{\oN^{2n}}$,
 its {\em complex moment
matrix} is the matrix $\fCM(y)$ indexed by $\oN^{2n}$
whose $(\alpha\alpha',\beta\beta')$th entry is
$y_{\alpha'+\beta,\alpha+\beta'}$,
for $(\alpha,\alpha'),(\beta,\beta')\in\oN^{2n}$. If $y$ satisfies
\begin{equation}\label{relher}
y_{\alpha'\alpha}=\overline{y_{\alpha\alpha'}} \ \
\text{ for }  (\alpha,\alpha')\in \oN^{2n},
\end{equation}
then  $\fCM(y)$ is a Hermitian matrix.
Let $M^\oC(y) $ denote the principal submatrix of $\fCM(y)$ indexed by
the subset $\{(0\alpha')\mid \alpha'\in \oN^n\}$; in other words,
one may think of $M^\oC(y)$ as being  indexed by $\oN^n$ with
$(\alpha',\beta')$th entry $y_{\alpha'\beta'}$; let us call $M^\oC(y) $ a {\em pruned complex moment matrix}.
These three types of matrices $M^\oK(y)$ ($\oK=\oR,\oC$) and $\fCM(y)$
will play a central role in our treatment.
It will be convenient to think of $M^\oK(y)$ as being indexed by
$\oT_{n}$ 
and of $\fCM(y) $ as being indexed by 
$$\ZZ_n:=
\{\oz^\alpha x^{\alpha'}\mid \alpha,\alpha'\in \oN^n\}\subseteq \oC[x,\oX].$$
Thus $\ZZ_n\sim \oT_{2n}$ and  we view $x$ as a complex variable in the complex case.
Recall that 
one says that `$f\in \PK$ lies in the kernel of $M^\oK(y)$' if $M^\oK(y) \ve{f}=0$.
Similarly, one may identify a polynomial $(x,\bar x)\mapsto f(x,\bar x)=\sum_{\alpha,\alpha'}f_{\alpha,\alpha'}
\bar x^\alpha x^{\alpha'}$ with its sequence of coefficients 
$\ve{f}=(f_{\alpha,\alpha'})_{\alpha,\alpha'}$ which allows us to say that 
`$f\in \oC[x,\oX]$ lies in $\Ker \fCM(y)$'
if $\fCM(y)\ve{f}=0$.

We also  need {\em truncated} moment matrices. For an integer $t\ge 0$, 
$M^\oK_t(y)$ denotes the 
principal submatrix of $M^\oK(y)$ indexed by $\oT_{n,t}$ 
and $\fCM_t(y)$
 denotes  the  principal submatrix of $\fCM(y)$
indexed by the set
$\ZZ_{n,t}:=\{\ox^\alpha x^{\alpha'}\mid \alpha,\alpha'\in \oN^n, \ |\alpha|+|\alpha'|\le t\}
\sim \oT_{2n,t}$.
Given $h\in \PR$, $h(x)=\sum_\beta h_\beta x^\beta$, 
 and $y\in \oR^{\oN^n}$, define $hy\in \oR^{\oN^n}$ by
$$hy:=M^\oR(y)\ve{h};\ \text{ that is, } (hy)_\alpha=\sum_\beta h_\beta y_{\alpha+\beta} \ \text{ for } \alpha\in \oN^n.$$
Similarly, given  $h(x,\oz)=\sum_{\beta,\beta'} h_{\beta\beta'} \oz^\beta x^{\beta'}\in \oC[x,\oX]$ and 
$y\in \oC^{\oN^{2n}}$,
define $hy\in \oC^{\oN^{2n}}$ by
$$hy:=\fCM(y)\ve{h}; \ \text{ that is, } 
(hy)_{\alpha\alpha'}=\sum_{\beta,\beta'}h_{\beta\beta'}
y_{\alpha'+\beta,\alpha+\beta'} 
\ \text{ for } \alpha,\alpha'\in\oN^n.$$
When $h\in \PC$ (i.e., $h_{\beta\beta'}=0$ if $\beta\ne 0$),
$M^\oC(y)\ve{h}$ is the projection  of $\fCM(y)\ve{h}$ onto the coordinates indexed by the pairs 
$(\alpha, \alpha')$ with 
$\alpha =0$.

\subsection{Measures and kernels of moment matrices}\label{secmeasure}
For a Hermitian matrix $A$, write $A\succeq 0$ if 
$A$ is positive semidefinite, i.e., if
$u^*Au\ge 0$ for all $u\in \oC^n$ (or $u\in \oR^n$ when $A$ is real valued).

\medskip\noindent
{\bf \textit{The real case.}}
For $v\in \oC^n$,  set $\zeta_v:=(v^\alpha)_{\alpha\in\oN^n}$ and
$\zeta_{t,v}:=(v^\alpha)_{\alpha\in\oN^n_t}$ for an integer $t\ge 0$.
Let $\mu$ be a positive measure on $\oR^n$ with finite support; 
 say, $\mu=\sum_{v\in W}\lambda_v\delta_v$ where $\lambda_v>0$ and $W\subseteq \oR^n$, $|W|<\infty$.
The {\em sequence of moments
of the measure $\mu$}  is the sequence $y^\mu\in \oR^{\oN^n}$ 
defined by
$(y^\mu)_\alpha:=\int x^\alpha d\mu =\sum_{v\in W}\lambda_v v^\alpha$ for $\alpha\in\oN^n$;
$(y^\mu)_0=\sum_{v\in W}\lambda_v$ is the total mass of the measure, equal to 1 if $\mu$ is a probability measure.
We have
$$
y^\mu=\sum_{v\in W}\lambda_v \zeta_v.$$
Moreover, $M^\oR(y^\mu)=\sum_{v\in W}\lambda_v \zeta_{v}\zeta_v^T\succeq 0$ and 
$$\Ker M^\oR(y^\mu)= \{f\in\PR\mid f(v)=0\ \forall v\in W\}=I(W),$$
$$\Ker M^\oR_t(y^\mu) = I(W)\cap \oR[x]_t $$
(which 
follows from the fact that $\ve{f}^TM^\oR_t(\zeta_{2t,v})\ve{f} =f(v)^2$
for $f\in \oR[x]_t$).
Given polynomials $h_1,\ldots,h_m\in\PR$, 
let $d_j,d$ be defined as in (\ref{reldj}) and, 
for $t\ge d$, set
\begin{equation}\label{setKR}
K_t^\oR:=\{y\in \oR^{\oN^n_{2t}} \mid y_0=1,\ M^\oR_t(y)\succeq 0, M^\oR_{t-d_j}(h_jy)=0\ (j=1,\ldots,m)\}.
\end{equation}
Then, $K^\oR_t$ is a convex set which contains the vectors $\zeta_{2t,v}$ for all $v\in V_\oR(I)$.
The following geometric observation, which indicates how the
real radical ideal  of $I$ relates to the kernel of moment matrices,
will play a central role in our approach.

\begin{lemma}\label{lemkerR}
Let $I=\langle h_1,\ldots,h_m\rangle \subseteq \PR$, $t\ge d$, and let $y\in K^\oR_t$ 
for which $\rank M^\oR_{t}(y)$ is maximum.  Then,
$\Ker M^\oR_{t}(y)\subseteq \Ker M^\oR_{t}(z)$ for all $z\in K^\oR_t$.
Moreover, $\Ker M^\oR_t(y)\subseteq I(V_\oR(I))$.
\end{lemma}

\begin{proof} 
Let $z\in  K^\oR_t$. Then, $y':=\frac{1}{2}(y+z)\in  K^\oR_t$ and
$\Ker  M^\oR_{t}(y')=\Ker  M^\oR_{t}(y)\cap \Ker  M^\oR_{t}(z) \subseteq 
\Ker  M^\oR_{t}(y)$. As $\rank  M^\oR_{t}(y)\ge \rank  M^\oR_{t}(y')$, 
equality $\Ker  M^\oR_{t}(y)\cap \Ker  M^\oR_{t}(z)=
\Ker  M^\oR_{t}(y)$ holds, which implies
$\Ker  M^\oR_{t}(y)\subseteq  \Ker  M^\oR_{t}(z)$.
  As $\zeta_{2t,v}\in K^\oR_t$ for all $v\in V_\oR(I)$, this implies $\Ker M^\oR_{t}(y) \subseteq \cap_{v\in V_\oR(I)} \Ker M^\oR_{t}(\zeta_{2t,v})$ which in turn is contained in $I(V_\oR(I))$.  $\hfill$ $\Box$ \end{proof}

\medskip

\medskip\noindent
{\bf \textit{The complex case.}}
Let $\mu$ be a positive measure on $\oC^n$ with finite support; that is,
$\mu=\sum_{v\in W}\lambda_v\delta_v$ where $\lambda_v>0$ and $W\subseteq \oC^n$, $|W|<\infty$.
One can now define the {\em doubly-indexed sequence of moments $y^\mu\in \oC^{\oN^{2n}}$
of the measure $\mu$}  by
$(y^\mu)_{\alpha\alpha'}:=\int \oz^\alpha x^{\alpha'}d\mu =\sum_{v\in W}\lambda_v
\ov^\alpha v^{\alpha'}$ for $\alpha,\alpha'\in\oN^n$.
Thus $y^\mu$ satisfies (\ref{relher}) and
$$y^\mu=\displaystyle\sum_{v\in W}\lambda_v \zeta_{\ov}\otimes \zeta_v. $$
Therefore, 
$\fCM(y^\mu)=\sum_{v\in W}\lambda_v  \zeta_{v}\otimes \zeta_{\ov}
\left(\zeta_{\ov}\otimes \zeta_{v}\right)^T\succeq 0$ and, in particular,
$M^\oC(y^\mu)=\sum_{v\in W} \lambda_v \zeta_{\ov}\zeta_v^T\succeq 0$.
Moreover, 
$$\Ker M^\oC(y^\mu)= \{f\in\PC\mid f(v)=0\ \forall v\in W\}=I(W),$$
$$\Ker \fCM(y^\mu) =\{f\in \oC[x,\oX]\mid f(v,\ov)=0\ \forall v\in W\}$$
(using the fact that 
$\ve{f}^*\fCM(\zeta_{\ov}\otimes \zeta_v)\ve{f}=|f(v,\ov)|^2$ for 
$f\in \oC[x,\oX]$). 
Given polynomials $h_1,\ldots,h_m\in\PC$, 
$t\ge d$, define the sets 
\begin{equation}\label{setKC}
\begin{array}{l}
K_t^\oC:=\{y\in\oC^{\oN^{2n}_{2t}} 
\mid y_0=1,\  (\ref{relher})
,\ M^\oC_t(y)\succeq 0,  \ \\ 
\phantom{K_t^\oC:=\{y\in\oC^{\oN^{2n}_{2t}}} M^\oC_{t-d_j}(h_jy)=0\ (j=1,\ldots,m)\},\\
\fCK_t:=\{y\in\oC^{\oN^{2n}_{2t}}\mid y_0=1, \ (\ref{relher})
,\  \fCM_t(y)\succeq 0,\ \\
\phantom{\fCK_t:=\{y\in\oC^{\oN^{2n}_{2t}}} \fCM_{t-d_j}(h_jy)=0 \ (j=1,\ldots,m)\}.
\end{array}
\end{equation}
Hence,  $\fCK_t\subseteq K_t^\oC$ are  both convex sets. The following analogue of Lemma~\ref{lemkerR}
holds in the complex case; we omit the proof.

\begin{lemma}\label{lemkerC}
Let $I=\langle h_1,\ldots,h_m\rangle \subseteq \PC$ and  $t\ge d$.
\begin{description}
\item[(i)] 
Let $y\in K^\oC_t$ 
for which $\rank M^\oC_{t}(y)$ is maximum. Then, 
$\Ker M^\oC_{t}(y)\subseteq \Ker M^\oC_{t}(z)$ for all $z\in K^\oC_t$.
Moreover, $\Ker M^\oC_t(y)\subseteq I(\VC(I))$.
\item[(ii)] Let $y\in \fCK_t$ for which $\rank \fCM_t(y)$ is maximum.
Then, 
    $\Ker \fCM_{t}(y)\subseteq \Ker \fCM_{t}(z)$ for all $z\in \fCK_t$.
Moreover, $\Ker M^\oC_t(y)\subseteq I(\VC(I))$.
\end{description}
\end{lemma}

\medskip\noindent
{\bf \textit{Link between the real and complex cases.}}
As shown e.g. in \cite{CF02} the complex moment problem in $\oC^n$ can be reduced to 
the real moment problem in $\oR^{2n}$.
Let us sketch the main idea. For $\alpha,\alpha'\in \oN^n$, define the polynomials
\begin{eqnarray*}
\Phi^{(\alpha\alpha')}(x,\ox)& :=&\left(\frac{x-\ox}{2i}\right)^\alpha 
\left(\frac{x+\ox}{2}\right)^{\alpha'}\,=\,
\sum_{\beta\beta'}\varphi^{(\alpha\alpha')}_{\beta\beta'} \ox^\beta
x^{\beta'}\\
\Psi^{(\alpha\alpha')}(u,v)&:=& (v-iu)^\alpha (v+iu)^{\alpha'}
\,=\,\sum_{\beta\beta'}\psi^{(\alpha\alpha')}_{\beta\beta'}u^\beta v^{\beta'}.
\end{eqnarray*}
The following can be easily verified: For all $x\in\oC^n$, and all $(u,v)\in\oR^n$,
\[{\overline{\Phi^{(\alpha\alpha')}(x,\ox)}} = \Phi^{(\alpha\alpha')}(x,\ox),\qquad
{\overline{\Psi^{(\alpha\alpha')}(u,v)}} = \Psi^{(\alpha'\alpha)}(u,v).
\]
In addition, $\Phi^{(\alpha\alpha' +\beta\beta')}=\Phi^{(\alpha\alpha')}\Phi^{(\beta\beta')}$, and
$\Psi^{(\alpha\alpha'+\beta\beta')} = \Psi^{(\alpha\alpha')}\Psi^{(\beta\beta')}$. Moreover, for every $\alpha,\alpha'\in\oN^n$,
\[\sum_{\beta\beta'}\psi^{(\alpha\alpha')}_{\beta\beta'}\Phi^{(\beta\beta')}(x,\ox)=\ox^\alpha x^{\alpha'};\quad
\sum_{\beta\beta'}\varphi^{(\alpha\alpha')}_{\beta\beta'}\Psi^{(\beta\beta')}(u,v)=u^\alpha v^{\alpha'}.\]
Next, given $y\in\oC^{\oN^{2n}}$, define the linear mapping $L_y:\oC[x,\ox]\rightarrow \oC$ by
$L_y(f)=\sum_{\beta\beta'}f_{\beta\beta'}y_{\beta\beta'}$ for $f\in \oC[x,\ox]$ with
$f(x,\ox)=\sum_{\beta\beta'} f_{\beta\beta'}\ox^\beta x^{\beta'}$, and 
the mapping $\varphi:\,\oC^{\oN^{2n}}\to\oC^{\oN^{2n}}$, $y\mapsto a:=\varphi(y)$ with:
\begin{equation}\label{link}
a_{\alpha\alpha'}=L_y(\Phi^{(\alpha\alpha')})=\sum_{\beta\beta'}\varphi^{(\alpha\alpha')}_{\beta\beta'} y_{\beta\beta'} \qquad  (\alpha\alpha'\in\oN^{2n}).
\end{equation}
Notice that $a\in \oR^{\oN^{2n}}$ if $y$ satisfies (\ref{relher}).
Conversely, given $a\in \oC^{\oN^{2n}}$, 
let $L_a: \oC[u,v]\rightarrow \oC$ be the linear mapping
\[g\:(=\sum_{\beta\beta'}g_{\beta\beta'}u^\beta v^{\beta'})\,\mapsto L_a(g)\,:=\,\sum_{\beta\beta'}g_{\beta\beta'}a_{\beta\beta'},\qquad g\in \oC[u,v],\] 
and  the linear mapping $\psi:\oC^{\oN^{2n}}\to \oC^{\oN^{2n}}$,
$a\mapsto y:=\psi(a)$ by
\[y_{\alpha\alpha'}\,=\, L_a(\Psi^{(\alpha\alpha')})\qquad
(\alpha\alpha'\in\oN^{2n}).\]
Notice that $y$ satisfies (\ref{relher}) whenever $a$ is real valued.
The mappings $\varphi$ and $\psi$ are inverse bijections between the set of sequences in
$\oC^{\oN^{2n}}$ satisfying (\ref{relher}) and $\oR^{\oN^{2n}}$.
Based on the above observations, we can now verify that
$$\fCM(y)\succeq 0 \Longleftrightarrow M^\oR(a)\succeq 0.$$
Assume first $\fCM(y)\succeq 0$ and let $f\in \oR^{\oN^{2n}}$ arbitrary. Then
\begin{eqnarray*}
f^TM^\oR(a)f&=&\sum_{\alpha\alpha',\beta\beta'}f_{\alpha\alpha'}f_{\beta\beta'}a_{\alpha\alpha'+\beta\beta'}
= L_y(\sum_{\alpha\alpha',\beta\beta'}f_{\alpha\alpha'}f_{\beta\beta'} \Phi^{(\alpha\alpha'+\beta\beta')})
\\
&=& L_y((\sum_{\alpha\alpha'}f_{\alpha\alpha'} \Phi^{(\alpha\alpha')})^2)
=\ve {g^*}\fCM(y)\ve g\ge 0,
\end{eqnarray*}
with $g(x,\ox):=\sum_{\alpha\alpha'}f_{\alpha\alpha'} \Phi^{(\alpha\alpha')}(x,\ox)$.
This shows that $M^\oR(a)\succeq 0$. Conversely, assume $M^\oR(a)\succeq 0$, and let $g\in \oC^{\oN^{2n}}$ arbitrary. Then
\begin{eqnarray*}
g^*\fCM(y)g&=&
\sum_{\alpha\alpha',\beta\beta'} \overline{g_{\alpha\alpha'}} g_{\beta\beta'} y_{\alpha'+\beta,\alpha+\beta'}\\
&=& L_a(\sum_{\alpha\alpha',\beta\beta'} \overline{g_{\alpha\alpha'}} g_{\beta\beta'}\, 
\Psi^{(\alpha'\alpha + \beta\beta')})
= L_a(\bar h h),
\end{eqnarray*}
with $h(u,v):= \sum_{\alpha\alpha'} g_{\alpha\alpha'} \Psi^{(\alpha\alpha')}(u,v)$.
Now $h=h_1+ih_2$ with $h_1,h_2\in\oR[u,v]$, and so
\begin{eqnarray*}
L_a(\bar h h)&=&L_a(h_1^2+h_2^2)\\
&=&\ve{ h_1}^TM^\oR(a)\ve{h_1}+ \ve{h_2}^T M^\oR(a)\ve{h_2}\ge 0,
\end{eqnarray*}
which shows that $\fCM(y)\succeq 0$.

Finally, $y$ is the sequence of moments of a measure on the set $W\subseteq \oC^n$ 
if and only if $a$ is the sequence of moments of a measure on the set
$\{(\text{Im}(v),\text{Re}(v))\mid v\in W\}\subseteq \oR^{2n}$.
(Use the fact that, if $y=\zeta_{\bar v }\otimes \zeta_v$, then $a=
\zeta_{(\text{Im}(v),\text{Re}(v))}$.)

\subsection{Flat extensions and finite rank  moment matrices}\label{secflat}

Given a Hermitian  matrix $A$ and a principal submatrix $B$ of $A$, one 
says that $A$ is a {\em flat extension} of $B$ if $\rank A=\rank B$;
then $A\succeq 0\Longleftrightarrow B\succeq 0$.
We begin with two fundamental results of Curto and Fialkow \cite{CF96} 
about {\em finite rank}  moment matrices, where this notion of flat extension plays a central role. See \cite{Lau05} for a
 short proof of Theorem \ref{theoCFmeasure}
and \cite{Lau2} for an exposition of Theorem~\ref{theoCFflat}.

\begin{theorem}\label{theoCFmeasure}
\begin{description}\item[(i)]
If $M^\oR(y)\succeq 0$ and $\rank M^\oR(y)<\infty$, then
$y=\sum_{v\in W}\lambda_v \zeta_v$ for some 
finite set $W\subseteq \oR^n$ and $\lambda_v>0$,  
$|W|=\rank M^\oR(y)$, and $\Ker M^\oR(y)= I(W)$.
\item[(ii)]
If $\fCM(y)\succeq 0$ and $\rank \fCM(y)<\infty$,
then $y=\sum_{v\in W} \lambda_v~\zeta_{\ov}\otimes\zeta_v$
for some
 finite set $W\subseteq \oC^n$ and $\lambda_v>0$,
 $|W|=\rank \fCM(y)$, and $\Ker M^\oC(y) = I(W)$. 
 \end{description}
 \end{theorem}

\begin{theorem}\label{theoCFflat}
\begin{description}\item[(i)]
If $M^\oR_t(y)\succeq 0$ and $\rank M^\oR_{t}(y)=\rank M^\oR_{t-1}(y)$, then 
$y$ can be extended in a unique way to $\tilde y\in \oR^{\oN^n}$ such that
$M^\oR(\tilde y)$ is a flat extension of $M^\oR_t(y)$ (and thus  $M^\oR(\tilde y) \succeq 0$).
\item[(ii)] 
If $\fCM_t(y)\succeq 0$ and $\rank \fCM_t(y)=\rank \fCM_{t-1}(y)$, then
$y$ can be extended in a unique way to $\ty\in \oC^{\oN^{2n}}$ such that 
$\fCM(\ty)$ is a flat extension of $\fCM_t(y)$ (and thus $\fCM(\ty)\succeq 0$).
\end{description}
\end{theorem}

The following lemma taken from  \cite{CF96} shows that the kernel 
of a truncated moment matrix enjoys  ideal-like
properties.

\begin{lemma}\label{lemideal}
\begin{description}
\item[(i)] Let $M^\oR_t(y)\succeq 0$, 
 $f,g\in \PR$, 
with $\deg(fg)\le t-1$. Then,
$M^\oR_t(y)\ve{f}=0\Longrightarrow M^\oR_t(y)\ve{fg}=0$. 
\item[(ii)] Let $\fCM_t(y)\succeq 0$, $f,g\in \oC[x,\oX]$,
with $\deg(fg)\le t-1$. Then, \\
$\fCM_t(y)\ve{f}=0\Longrightarrow \fCM_t(y)\ve{fg}=0$.
\end{description}
\end{lemma}

\begin{proof} Set $h:=fg$. 
(i) As $\deg(h)\le t-1$ and $M_t^\oR(y)\succeq 0$, it suffices to show
$M^\oR_{t-1}(y)\ve{h}=0$.
Moreover it suffices to show the result for $g=x_i$; in this latter case one can verify
that, for $\alpha\in\oN^n_{t-1}$,
$(M_{t-1}^\oR(y)\ve{h})_\alpha= (M_t^\oR(y)\ve{f})_{\alpha+e_i}=0$.\\
(ii) Similarly assume $g=x_i$ or $\overline x_i$.
For $\alpha\alpha'\in \ZZ_{n,t-1}$, 
$(\fCM_{t-1}(y)\ve{h})_{\alpha\alpha'}$ is equal to $(\fCM_t(y)\ve{f})_{\alpha+e_i\ \alpha'}$
if $g=x_i$ and to 
$(\fCM_t(y)\ve{f})_{\alpha\  \alpha'+e_i}$ if $g=\overline x_i$, thus to 0 in both cases.
$\hfill\Box$
\end{proof}

\begin{proposition}\label{lemmomR} 
$\Ker M^\oR(y)$ is an ideal in $\PR$, which is 
real radical if $M^\oR(y)\succeq 0$. 
Assume $M^\oR(y)\succeq 0$ and
$\rank M^\oR(y)=\rank M^\oR_{t-1}(y)$ for some integer $t\ge 1$.
Then, $\Ker M^\oR(y)=\langle \Ker M^\oR_{t}(y)\rangle$ and,
for $\BB\subseteq \oT_{n}$,
$\BB$  indexes a (maximum) nonsingular principal submatrix of $M^\oR(y)$
if and only if $\BB$ is a (maximum) linearly independent 
subset of $\PR/\Ker M^\oR(y)$.
\end{proposition}

\begin{proof} 
We use the (easy to verify) identity: 
$$\ve{h}^TM^\oR(y) {\ve{pq}}= \ve{hq}^T M^\oR(y)\ve{p}$$
for  $p,q,h\in \PR$.
If $p\in \Ker M^\oR(y)$, $q\in \PR$, then $\ve{h}^TM^\oR(y) \ve{pq}= \ve{hq}^T M^\oR(y) \ve{p}=0$
for all $h\in \PR$, which implies  $M^\oR(y) \ve{pq}=0$ and thus $pq\in \Ker M^\oR(y)$.
This shows that 
$\Ker M^\oR(y)$ is an ideal.
Assume now $M^\oR(y)\succeq 0$; we show that $\Ker M^\oR(y)$ is  real radical. In  view of Lemma \ref{lemradical}, it suffices to show that
if $\sum_{i=1}^k p_i^2\in \Ker M^\oR(y)$ for some $p_i\in\PR$, then 
$p_i\in 
\Ker M^\oR(y)$. 
Indeed,  $0=\ve{1}^T M^\oR(y) \ve{\sum_{i=1}^k p_i^2}= 
 \sum_{i=1}^k \ve{p_i}^T M^\oR(y) \ve{p_i}$ implies $\ve{p_i}^T M^\oR(y) \ve{p_i}=0$ and thus 
$p_i\in \Ker M^\oR(y)$ for all $i$.

Assume $\rank M^\oR(y)=\rank M^\oR_{t-1}(y)=:r$
and set $J:=\langle \Ker M^\oR_{t}(y)\rangle$. Obviously, $J\subseteq \Ker M^\oR(y)$;
we  show equality.
For this, let $\BB\subseteq \oT_{n,t-1}$ 
index an $r\times r$  nonsingular principal submatrix of $M^\oR(y)$. 
We show that, for all $\alpha\in\oN^n$,
$x^\alpha \in \Span_\oR(\BB)+J$, using induction on $|\alpha|$. This holds for $|\alpha|\le t$
by the definition of $\BB$.
Assume $|\alpha|\ge t+1$ and write $x^\alpha=x_ix^\delta$.
By the induction assumption,
$x^\delta = \sum_{x^\beta\in\BB} c_\beta x^\beta + q$ where $q\in J$, $c_\beta\in\oR$.
Thus, $x^\alpha = \sum_{x^\beta\in\BB} c_\beta x_ix^\beta  +x_iq$.
Here, $x_iq\in J$ and $x_ix^\beta \in \Span_\oR(\BB)+J$ since $\deg(x_ix^\beta)\le t$, which implies $x^\alpha \in \Span_\oR(\BB)+J$.
Thus we have shown that $\PR=\Span_\oR(\BB)+J$. 
As $\Ker M^\oR(y)\cap \Span_\oR(\BB)=\{0\}$, this implies easily 
that $\Ker M^\oR(y)=J$.

For $\BB\subseteq \oT_n$, it is obvious that $\BB$ indexes a nonsingular submatrix
of $M^\oR(y)$ if and only if $\BB$ is linearly independent in
$\PR/\Ker M^\oR(y)$. The last statement of the lemma now follows 
since $\dim\PR/\Ker M^\oR(y) = r$ 
(as $\Ker M^\oR(y)$ is radical and using the identity
$ |\VC(\Ker M^\oR(y))| = \rank M^\oR(y)$ from Theorem \ref{theoCFmeasure}).
$\hfill$ $\Box$ 
\end{proof}

\begin{proposition}\label{lemmomC}
$\Ker \fCM(y)$ is an ideal in $\oC[x,\oX]$. 
If $\fCM(y)\succeq~0$,
then $\Ker \fCM(y)$ is a radical ideal in $\oC[x,\oX]$ and thus $\Ker M^\oC(y)$ is a radical ideal in $\PC$.
Assume, moreover, $\rank M^\oC(y)=\rank M^\oC_{t-1}(y)$ for some integer $t\ge 1$.
Then, $\Ker M^\oC(y)=\langle \Ker M^\oC_t(y)\rangle$ and, 
for  
$\BB\subseteq \oT_{n}$, $\BB$  indexes a (maximum) nonsingular principal submatrix
of $M^\oC(y)$ if and only if $\BB$ is  a (maximum) linearly independent
subset of $\oC[x]/\Ker M^\oC(y)$.
\end{proposition}

\begin{proof}
As in the real case, we use 
the following (easy to verify)
identities: For
$h,p,q\in\oC[x,\bar x]$,
$$\begin{array}{l}
\ve{h}^*\fCM(y)\ve{pq} = \ve{h\overline p}^*\fCM(y)\ve{q},\\
\ve{p}^*\fCM(y)\ve{p}=\ve{1}^*\fCM(y)\ve{p\overline{p}},\\
\ve{p^2}^*\fCM(y) \ve{p^2} = \ve{p\overline{p}}^*\fCM(y)\ve{p\overline{p}},
\end{array}$$
where $\overline p\in \oC[x,\ox]$ is defined as $\overline p(x,\ox):=\overline {p(x,\ox)}$. 
This implies directly that $\Ker \fCM(y)$ is an ideal.
Assume now $\fCM(y)\succeq 0$;
we show that ${\rm Ker}M^{2\oC}(y)$ is radical. 
In view of Lemma \ref{lemradical}, this follows from the following fact:
$$\begin{array}{l}
p^2\in{\rm Ker}M^{2\oC}(y) \\
\Longrightarrow 0=\ve{p^2}^*\fCM(y)\ve{p^2} =\ve{p\overline{p}}^*\fCM(y)\ve{p\overline{p}} \\
\Longrightarrow p\overline{p}\in \Ker \fCM(y) \\
\Longrightarrow 0= \ve{1}^*\fCM(y)\ve{p\overline{p}}=\ve{p}^*\fCM(y)\ve{p} \\
\Longrightarrow p\in  \Ker \fCM(y).
\end{array}$$
The proof for the last statements of the proposition is identical to the proof of the corresponding statements in Proposition \ref{lemmomR} for the real case.
$\hfill\Box$ \end{proof}

\medskip
Without the assumption $\fCM(y)\succeq 0$, 
$\Ker M^\oC(y)$ is not necessarily an ideal in $\PC$. Indeed, for the sequence 
$y\in \oC^{\oN^{2n}}$ defined by $y_{\alpha\alpha'}:=0$ if $\alpha=0$ or
$\alpha'=0$, and
$y_{\alpha\alpha'}:=1$ otherwise,  $M^\oC(y)\succeq 0$,
$\fCM(y)\not\succeq 0$ and
$\Ker M^\oC(y)$  is not an ideal (e.g.,  $1\in \Ker M^\oC(y)$
while any nonconstant monomial does not lie in  $\Ker M^\oC(y)$).
We mention for further reference the following corollary,
 and we conclude the section with the proof of Proposition \ref{prop0} and a lemma about the `(real) radical'-like property of the kernel of a positive semidefinite truncated moment matrix.

\begin{corollary}\label{corker}
\begin{description}
\item[(i)]
Assume $M^\oR_s(y)\succeq 0$ and
$\rank M_s^\oR(y)=\rank M_{s-1}^\oR(y)=:r$.
Then, $J:=\langle \Ker  M_s^\oR(y)\rangle$ is real radical and zero-dimensional, $\dim \PR/J=r$,
$J \cap \oR[x]_s
=\Ker M_s^\oR(y)$ and,
for $\BB\subseteq \oT_{n,s}$, $\BB$ indexes a (maximum) nonsingular principal submatrix
of $M^\oR_s(y)$ $\Longleftrightarrow$ 
$\BB$ is (maximum) linear independent in $\PR/J$.
\item[(ii)] Assume $\fCM_s(y)\succeq 0$ and  $\rank \fCM_s(y) =
\rank M^\oC_{s-1}(y)=:r$.
Then, $J:=\langle \Ker  M_s^\oC(y)\rangle$ is radical, $\dim \PC/J=r$, $J\cap\oC[x]_s=\Ker M^\oC_s(y)$ and, for
$\BB\subseteq \oT_{n,s}$, $\BB$ indexes a (maximum) nonsingular principal submatrix
of $M^\oC_s(y)$ $\Longleftrightarrow$
$\BB$ is (maximum) linearly independent in $\PC/J$.
\end{description}
\end{corollary}

\begin{proof} We prove only (i). By Theorem \ref{theoCFflat}, $y$ has an extension 
$\tilde y\in \oR^{\oN^n}$ such that $M^\oR(\tilde y)$ is a flat extension
of $M^\oR_t(y)$. By Proposition \ref{lemmomR}, the ideal 
$\Ker M^\oR(\tilde y)=\langle \Ker M_s(y)\rangle =:J$ is real radical and zero-dimensional,
$\dim \PR/J=r$, and 
$J\cap\oR[x]_s =\Ker M^\oR(\tilde y)\cap \oR[x]_s=\Ker M^\oR_s(\tilde y)=
\Ker M^\oR_s(y)$. 
$\hfill\Box$
\end{proof}

\medskip\noindent
{\sc  Proof of Proposition \ref{prop0}.}
By Proposition \ref{lemmomR}, $J:=\Ker M^\oR(y)$ is a real radical ideal, since $M^\oR(y)\succeq 0$.
As $0=M^\oR(h_jy)=M^\oR(y)\ve{h_j}$ for all $j$, we have $I\subseteq J$, which implies that 
$V_\oR(J)\subseteq V_\oR(I)$ is finite.
As $J$ is real radical, we deduce that $\VC(J)=V_\oR(J)\subseteq \oR^n$. Hence $J$ is zero-dimensional and
$I(V_\oR(I))\subseteq J$ since $\VC(J)\subseteq V_\oR(I)$.
Set $r:=\dim \PR/J = |\VC(J)|\le |V_\oR(I)|.$
Let $\BB\subseteq \oT_n$ be a linear basis of 
$\PR/J$, $|\BB|=r$. Then the columns of $M^\oR(y)$ indexed by $\BB$ form a basis of the column space of $M^\oR(y)$ and thus $\rank M^\oR(y)=r$.
Moreover,
$r=|V_\oR(I)|$ if and only if $\VC(J)=V_\oR(I)$ which in turn is equivalent to 
$J=I(V_\oR(I))$. Now, this maximum rank $|V_\oR(I)|$ is reached 
by the sequence $y:=y^\mu=\sum_{v\in V_\oR(I)} \lambda_v \zeta_v$ 
with $ \lambda_v >0$ which indeed satisfies (\ref{eq0}).
$\hfill\Box$

\medskip
Finally, it is useful to observe that the kernel of a positive semidefinite truncated moment matrix enjoys 
the following `(real) radical'-like property.
We omit the proof whose details are straightforward.

\begin{lemma}\label{lem0}
\begin{description}
\item[(i)] 
Assume $M^\oR_t(y)\succeq 0$ and let  $p,q_j\in \PR$, $f:= p^{2m}+\sum_jq_j^2$ with 
 $m\in \oN$, $m\ge 1$. Then,
 $f\in \Ker M_t^\oR(y)\Longrightarrow 
p\in \Ker M_t^\oR(y)$.
\item[(ii)] Assume $\fCM_t(y)\succeq 0$ and let $p\in \PC$, $m\in \oN$, $m\ge 1$.
Then, $p^m\in \Ker \fCM_t(y)\Longrightarrow p\in \Ker \fCM_t(y)$.
\end{description}
\end{lemma}

\section{A semidefinite characterization of the (real) radical ideal via moment matrices}\label{secfind}

In this section we present a semidefinite characterization of the real radical 
ideal $I(V_\oR(I))$ of an ideal $I\subseteq\oR[x]$, as well as a numerical
algorithm for computing a set of generators. It turns out that the method also 
applies to the radical ideal $I(\VC(I))$. Our strategy is to obtain $I(V_\oK(I))$ ($\oK=\oR$ or $\oC$) as the ideal generated by the kernel of some suitable moment matrix $M^\oK_t(y)$ where $y\in K^\oK_t$. Sections \ref{secfind1}-\ref{secfind3} contain some results ensuring that the moment matrix $M^\oK_t(y)$ has the desirable properties for achieving this task and Section \ref{subsec::algorithm} describes our algorithm. 

\subsection{Weakest set of conditions}\label{secfind1}
Throughout, $I=\langle h_1,\ldots,h_m\rangle$ is an ideal in $\PK$ for which we want to find 
the radical ideal $I(V_\oK(I))$, $\oK=\oR$ or $\oC$. Recall the definition of $d$ in (\ref{reldj}).

\begin{proposition}\label{propborder}
Let $t\ge d$, $1\le s\le t$, $y\in K_t^\oK$ for which $\rank M^\oK_t(y)$ is maximum,
 and 
$\BB\subseteq \oT_{n,s-1}$ index a maximum nonsingular principal 
submatrix of $M^\oK_{s-1}(y)$ with  border $\pBB$ defined as in (\ref{relborder}).
Assume (i)-(iii) below hold:
\begin{description}
\item[(i)] $\BB$ is an order ideal.
\item[(ii)] The principal submatrix of $M_{s}^\oK(y)$ indexed by $\BB\cup \pBB$ 
has the same rank as $M_{s-1}^\oK(y)$; that is,
 with $\BB:=\{b_1,\ldots,b_N\}$ and
$\pBB:=\{c_1,\ldots,c_H\}$, 
there exists  a polynomial $g_j\in \Ker M_{s}^\oK(y)$ of the form 
$g_j(x)=c_j-\sum_{i=1}^Na_{ij}b_i$
(i.e., $G:=\{g_1,\ldots,g_H\}$ is a $\BB$-border prebasis).
\item[(iii)] The formal multiplication matrices $\XX_1,\ldots,\XX_n$ defined from 
$G$ commute pairwise. 
\end{description}
Then $G$ is a border basis of $J:=\langle G\rangle\subseteq I(V_\oK(I))$,
$\BB$ is a linear basis of $\PK/J$,
and one can 
extract (using the formal multiplication matrices)
the  set $W:=\VC(J)$ which satisfies
 $V_\oK(I)\subseteq W$ and $|W|\le \rank M_{s-1}^\oK(y)$. 
Moreover, if $|V_\oK(I)|=|W|=\rank  M_{s-1}^\oK(y)$, then 
$V_\oK(I)=W$, $I(V_\oK(I))=J$, and $G$ is a $\BB$-border basis of 
$I(V_\oK(I))$.
\end{proposition}

\begin{proof}
Theorem \ref{theoborder} gives directly 
that $G$ is a border basis of 
the ideal $J=\langle G\rangle$ and  that
 $\BB$ is a linear basis of $\oR[x]/J$. Moreover,
the matrices $\XX_1,\ldots,\XX_n$ coincide with the multiplication matrices 
in $\oR[x]/J$ w.r.t. the basis $\BB$. Thus one can compute the set $W=\VC(J)$ from their eigenvectors.
By construction, $J\subseteq \langle \Ker M_{s}^\oK(y)\rangle
\subseteq \langle \Ker M_t^\oK(y)\rangle \subseteq 
I(V_\oK(I))$, where the last inclusion follows from  Lemmas \ref{lemkerR} and \ref{lemkerC} (i). 
This implies  $V_\oK(I)\subseteq \VC(J)=W$.
Moreover, $|W|\le \dim \PK/J=|\BB|=\rank M_{s-1}^\oK(y)$.

If $|V_\oK(I)|=|W|=\rank  M_{s-1}^\oK(y)$, then
$W=V_\oK(I)$ and $J$ is radical since 
$ \dim \PK/J = |\VC(J)|$, which implies 
 $I(V_\oK(I))=I(W)=I(\VC(J))=J$. 
$\hfill\Box$\end{proof}

\medskip
In the next two subsections we 
give  simple rank conditions (\ref{flati}),  (\ref{flatiC}), which  ensure that the 
conditions of Proposition \ref{propborder} hold. (See Remark \ref{comments} for details.)
For the sake of clarity, we treat the real and complex cases separately.

\subsection{Characterizing and computing the real radical $I(V_\oR(I))$}\label{secfind2}

Assume $I=\langle h_1,\ldots,h_m\rangle$ is an ideal in $\PR$, with $h_1,\ldots,h_m\in \PR$.

\begin{proposition}\label{propborderR}
Let $t\ge d$ and $y\in K^\oR_t$  for which $\rank M^\oR_{t}(y)$ is maximum. 
If,  for some $1\le s\le t$,
\begin{equation}\label{flati}
\rank M^\oR_{s}(y)=\rank M^\oR_{s-1}(y)=:r
\end{equation}
then $I(V_\oR(I))\supseteq \langle \Ker M^\oR_s(y)\rangle=:J$.
One can compute the  set $W:=\VC(J)$ which satisfies 
$ V_\oR(I)\subseteq W$  and $|W|=r$.
Moreover, if $V_\oR(I)=W$ then $I(V_\oR(I))=J$.
\end{proposition}

\begin{proof}
By Lemma \ref{lemkerR}, $\Ker M^\oR_s(y)\subseteq 
\Ker M^\oR_t(y) \subseteq I(V_\oR(I))$, implying 
$J\subseteq I(V_\oR(I))$ and thus $V_\oR(I)\subseteq W$.
By Corollary \ref{corker} (i), as $J$ is radical, 
$|W|=\dim\PR/J=r$, and
 $V_\oR(I)=W$ implies $I(V_\oR(I))=I(W)=J$.
$\hfill\Box$ \end{proof}

\medskip
One can verify (using Lemma \ref{lemideal}) that, if (\ref{flati}) holds for some 
$s\le t-2$, then it also holds for $s=t-1$. Hence it suffices to check
 whether (\ref{flati}) holds for $s=t-1$ or $t$.
In Lemma \ref{lemobs} below, we
observe that, if assumption (ii) in Proposition~\ref{propborder}
holds for $s\le t-1$,  
then in fact (\ref{flati}) holds and thus Proposition \ref{propborderR} applies.  
However, it may be that Proposition \ref{propborder} applies to the case $s=t$ while (\ref{flati}) does not hold; see Example~\ref{gauss} (for relaxation order $t=2$) for such an instance. 
In Remark \ref{comments} below we see that the converse of the next lemma holds.

\begin{lemma}
\label{lemobs}
In Proposition \ref{propborder}, if assumption (ii) holds for $s\le t-1$ then $\rank M_{s-1}^\oK(y)=\rank M_{s}^\oK(y)$.
\end{lemma}

\begin{proof} We have to show that $\oT_{n,s}\subseteq  
\Span_\oK(\BB)+\Ker M_s(y)$.
 By the definition
of $\BB$, any  $x^\alpha\in \oT_{n,s-1}$ lies in  $\Span_\oK(\BB)+\Ker M_{s-1}(y)$.
 For $x^\alpha\in\oT_{n,s}$,
write $x^\alpha = x_1x^\delta$ where $x^\delta\in\oT_{n,s-1}$.
 Thus $x^\delta=\sum_{x^\beta\in \BB}
c_\beta x^\beta +p$ where $p\in \Ker M_{s-1}^\oK(y)$, $c_\beta\in\oK$.
Therefore,
$x^\alpha =\sum_{x^\beta\in\BB} c_\beta x_1x^\beta +x_1p$.
As $x_1x^\beta\in\BB\cup\pBB$,  assumption (ii) implies 
that $x_1x^\beta\in\Span_\oK(\BB)+\Ker M_s^\oK(y)$.
As $\deg(x_1p)\le s\le t-1$, it follows from Lemma \ref{lemideal}
that $x_1p\in \Ker M_s^\oK(y)$. Therefore,
$x^\alpha \in \Span_\oK(\BB)+\Ker M_s(y)$.
$\hfill\Box$ \end{proof}

\bigskip
By strengthening the rank condition (\ref{flati}), one can show that 
$J=I(V_\oR(I))$, i.e., the desired real radical ideal is found.

\begin{proposition}\label{propcritR}
Let $t\ge d$ and $y\in K^\oR_t$ for which $\rank M^\oR_{t}(y)$ is maximum. 
Assume that, either
 (\ref{flati}) holds
 for some $2d\le s\le t$, or 
\begin{equation}\label{flatii}
\rank M^\oR_{s}(y)=\rank M^\oR_{s-d}(y)
\end{equation}
for some $d\le s\le t$.
Then $I(V_\oR(I))=\langle \Ker M^\oR_s(y)\rangle$ (and one
 can find $V_\oR(I)$). Moreover, $\rank M^\oR_s(y)=\vert V_\oR(I)\vert$.
\end{proposition}

\begin{proof} 
In view of Proposition \ref{propborderR}, there 
remains only to show the  inclusion $I(V_\oR(I))\subseteq J:=\langle \Ker M^\oR_s(y)\rangle$ or, equivalently (since $J$ is radical),
 $W:=\VC(J)\subseteq V_\oR(I)$. We already know that $W\subseteq \oR^n$ since $J$ is real radical and zero-dimensional
(see Corollary \ref{corker}). We now show that $W\subseteq \VC(I)$.
Assume first that (\ref{flati}) holds for $s\ge 2d$.
As $M^\oR_{t-d_j}(h_jy)=0$, we have $(h_jy)_\alpha=0$ for all $|\alpha|\le 2t-2d_j$
and thus for all $|\alpha|\le 2d_j$. Hence,
$M^\oR_{2d_j}(y)\ve{h_j}=0$ and thus $h_j\in \Ker M^\oR_s(y)$. Therefore,
$I\subseteq J$, giving $W\subseteq \VC(I)$.

Assume now that (\ref{flatii}) holds for some $d\le s\le t$.
Let $p_v$ ($v\in W$) be interpolation polynomials, i.e., $p_v(w)=\delta_{v,w}$ for $v,w\in W$. As observed in 
\cite[Lemma 27]{Lau2}, one can assume that $\deg(p_v)\le s-d$. (Indeed,
let $\BB\subseteq \oT_{n,s-d}$ index a maximum nonsingular submatrix
of $M^\oR_s(y)$; then $\BB$ is a basis of 
$\PR/J$ by Corollary \ref{corker} and one can replace $p_v$ by its residue modulo $J$ w.r.t. $\BB$.)
>From Theorems \ref{theoCFmeasure}, \ref{theoCFflat}, 
we know that $(y_\alpha)_{\alpha\in \oT_{n,2s}}
=\sum_{v\in W}\lambda_v \zeta_{2s,v}$ where $\lambda_v>0$.
Hence, $0=\ve{p_v}^TM^\oR_{s-d_j}(h_jy)\ve{p_v}=h_j(v)\lambda_v$ implies
$h_j(v)=0$ for all $j$ and thus
$v\in \VC(I)$, which shows $W\subseteq \VC(I)$.
$\hfill$ $\Box$ 
\end{proof} 

\medskip
We now formulate an analogous result for the ideal $I(V_\oR(I)\cap S)$, 
where 
$S:=\{x\in \oR^n\mid h_{m+1}(x)\ge 0,\ldots,h_{m+k}(x)\ge 0\}$ is a semialgebraic set, 
with $h_{m+1},\ldots,h_{m+k}\in \PR$.
For this define 
 the set 
\begin{equation}\label{setKtS}
K^\oR_{t,S}:=K^\oR_t\cap\{y\mid M_{t-d_j}(h_jy)\succeq 0\ (j=m+1,\ldots,m+k)\}
\end{equation}
for $t\ge d:=\displaystyle \max_{j=1,\ldots,m+k}d_j$.

\begin{proposition}\label{propborderS}
Let $t\ge d$ and $y\in K^\oR_{t,S}$ for which $\rank M^\oR_t(y)$ is maximum. 
\begin{description}
\item[(i)] 
Assume (\ref{flati}) holds for some $1\le s\le t$. Then
$J:=\langle \Ker M_s^\oR(y)\rangle \subseteq I(V_\oR(I)\cap S)$ and
$W:=\VC(J)\supseteq V_\oR(I)\cap S$ with $|W|=\rank M^\oR_s(y)$; moreover, 
$J= I(V_\oR(I)\cap S)$  if $W= V_\oR(I)\cap S$.
\item[(ii)] 
Assume (\ref{flatii}) holds for some $d\le s\le t$. Then
$I(V_\oR(I)\cap S)=\langle \Ker M_s^\oR(y)\rangle.$
\end{description}
\end{proposition}

\begin{proof}
(i) The inclusion $\Ker M_t^\oR(y)\subseteq  I(V_\oR(I)\cap S)$ 
follows from the maximality of the rank of $ M_t^\oR(y)$ and the fact that
$\zeta_{2t,v}\in K^\oR_{t,S}$ for all $v\in V_\oR(I)\cap S$. This gives
$J\subseteq  I(V_\oR(I)\cap S)$  and thus $W\supseteq V_\oR(I)\cap S$.
Equality $W=V_\oR(I)\cap S$ implies $I(V_\oR(I)\cap S)=I(W)=J$ (as $J$ radical).
This concludes the proof of (i). The proof for (ii) is analogous to that of 
the corresponding statement in Proposition \ref{propcritR}.
$\hfill\Box$
\end{proof}

\medskip
To conclude we  show that, when $V_\oR(I)$ is finite, then
condition (\ref{flatii}) is satisfied for $t$ large enough.
That is, the conclusion of Propositions \ref{propcritR}, \ref{propborderS} holds: the real radical
ideal $I(V_\oR(I))$ or $I(V_\oR(I)\cap S)=\langle \Ker M^\oR_s(y)\rangle$ is found.

\begin{proposition} \label{theorealidealR}
Assume $|V_\oR(I)|<\infty$.
\begin{description}
\item[(i)] If $V_\oR(I)=\emptyset$ then $K_{t,S}^\oR=\emptyset$ for $t$ large enough.
\item[(ii)] If $V_\oR(I)\ne \emptyset$ then, for $t$ large enough,
there exists $d\le s\le t$ such that
$\rank M_s^\oR(y)=\rank M_{s-d}^\oR(y)$ for all $y\in K_{t,S}^\oR$.
\end{description}
\end{proposition}

\begin{proof}
Assume $t\ge 2d$ and let $y\in K^\oR_{t,S}$. Then, as observed in the proof of Proposition \ref{propcritR},
 $h_1,\ldots,h_m\in \Ker M^\oR_t(y)$. We first show that, for $t$ large enough,
 $\Ker M^\oR_t(y)$ also contains a
given  basis of the ideal $I(V_\oR(I))$.

\begin{claim}\label{claim1}
 Let $\{g_1,\ldots,g_k\}$ be a basis of the ideal $I(V_\oR(I))$.
There exists $t_0\in\oN$ such that $g_1,\ldots,g_k\in\Ker M^\oR_t(y)$ for all
$t\ge t_0$.
\end{claim}

\begin{proof}
Let $l \in\{1,\ldots,k\}$. 
By the Real Nullstellensatz,
 there exist  $m_l\in\oN$, $m_l\ge 1$  and polynomials
$\sigma_l$, $u^{(l)}_j$ ($j\le m$) for which
$g_l^{2m_l} +\sigma_l=\sum_{j=1}^m u^{(l)}_j h_j$ and $\sigma_l$ is a sum of squares.
Set $t_0:= 1+\max_{l\le k, j\le m}(2d, \deg(g_l^{2m_l}), \deg(\sigma_l),$ $  \deg(u^{(l)}_jh_j))$
 and let $t\ge t_0$. As $\deg(u^{(l)}_jh_j)\le t-1$ and  $h_j\in \Ker M^\oR_t(y)$,   
then $u^{(l)}_jh_j\in \Ker M^\oR_t(y)$ by Lemma \ref{lemideal}.
Hence, $g_l^{2m_l}+\sigma_l\in \Ker M^\oR_t(y)$ which, using
Lemma \ref{lem0} (i), implies  $g_l\in \Ker M^\oR_t(y)$.
$\hfill\Box$
\end{proof}

\medskip
If $V_\oR(I)=\emptyset$, then $\{1\}$ is a basis of $I(V_\oR(I))=\PR$. 
Hence $1\in\Ker M_t^\oR(y)$, implying $y_0=0$, which contradicts
 the fact that $y_0=1$ for $y\in K_{t,S}^\oR$. 
Therefore, $K_{t,S}^\oR=\emptyset$ for $t\ge t_0$, which shows (i).

\medskip
Assume now $V_\oR(I)\ne\emptyset$. Let $\{g_1,\ldots,g_k\}$ be a Gr\"obner basis 
of $I(V_\oR(I))$ for a graded monomial ordering ordering and let $\BB$ be the corresponding 
set of standard monomials. Thus $\BB$ is a basis of $\PR/I(V_\oR(I))$; set 
$d_\BB:=\max_{b\in \BB} \deg(b)$ (which is well defined as $\BB\ne \emptyset$).
We can write any monomial as 
$x^\alpha = r^{(\alpha)}+ \sum_{l=1}^k p_l^{(\alpha)}g_l$, where
$r^{(\alpha)}\in \Span_\oR(\BB)$, $p_l^{(\alpha)}\in\PR$ and
 $\deg(p_l^{(\alpha)}g_l)\le \deg(x^\alpha)$.  
Set $t_1:=\max(d_\BB+d,t_0)$ and  let $t\ge t_1+1$. 
Consider $\alpha\in \oT_{n,t_1}$.
As $\deg(p_l^{(\alpha)}g_l)\le t_1\le t-1$ and $g_l\in\Ker M_t^\oR(y)$,
we have  $p_l^{(\alpha)}g_l \in \Ker M_t(y)$ and thus $x^\alpha-r^{(\alpha)}\in
 \Ker M_t(y)$. As $\deg(r^{(\alpha)})\le d_\BB\le t_1-d$, this shows that the $\alpha$th column of
$M_t^\oR(y)$ is a linear combination of columns indexed by
$\oT_{n,t_1}$. Therefore, $\rank M_{t_1}^\oR(y)=
\rank M_{t_1-d}^\oR(y)$, thus proving (ii).
$\hfill\Box$
\end{proof}

\medskip
\begin{remark}\label{remdetect}  {\bf Detecting existence of real solutions.}
Hence one can detect the existence of real solutions via the following criterion:
\begin{equation}\label{testsolution}
V_\oR(I)=\emptyset \Longleftrightarrow K_t^\oR=\emptyset \ \ \text{ for some } t.
\end{equation}
(The `only if' part follows directly 
from Proposition \ref{theorealidealR} (i), while the `if part' follows from the fact that $\zeta_{2t,v}\in K_{t}^\oR$ for any $v\in V_\oR(I)$.)
Note moreover that, when  $V_\oR(I)=\emptyset$, none of  the flat conditions (\ref{flatii}), or (\ref{flati}) with $s\ge 2d$, can hold;
indeed, under either of these two conditions, one would have
$|V_\oR(I)|=\rank M_s^\oR(y)\ge 1$ by Proposition \ref{propcritR}.
Consider as an illustration the following small example: $I=\langle h:=x_1^2+1\rangle\subseteq \oR[x_1]$ with
$V_\oR(I)=\emptyset$.
For $t\ge 1$, if $y\in K_t$ then $M_1(y)\succeq 0$ implies $y_{2e_1}\ge 0$, while $(hy)_0=0$ gives 
$y_{2e_1}+1=0$, yielding a contradiction. Hence,  $K_t=\emptyset$ for any $t\ge 1$.
\end{remark}

\medskip
\begin{remark}
Proposition \ref{theorealidealR} remains valid under the weaker
assumption $|V_\oR(I)\cap S|<\infty$ if, in the definition of the set
$K_{t,S}^\oR$ in (\ref{setKtS}), we add the constraints 
$M_{t-d_e}(p_ey)\succeq 0$ for $e\in \{0,1\}^k$, after setting
$p_e:= \prod_{i=1}^k h_{m+i}^{e_i}.$
The proof is analogous,
except we now prove in Claim \ref{claim1} that $\Ker M_t^\oR(y)$ contains a given
    basis of the ideal $I(V_\oR(I)\cap S)$. To show this, 
 instead of the Real Nullstellensatz, we now 
 use the Positivstellensatz (see Stengle \cite{St74}) which in our case can be formulated in the following way:
For $g\in\PR$, $g\in I(V_\oR(I)\cap S)$ if and only if 
$-g^{2r}=\sum_{j=1}^mu_jh_j +\sum_{e\in \{0,1\}^k}\sigma_e p_e$
for some $r\in \oN\setminus\{0\}$, $u_j,\sigma_e\in \PR$, with $\sigma_e$  s.o.s.
\end{remark}

\subsection{Characterizing and computing the radical $I(\VC(I))$}\label{secfind3}
 
Using complex moment matrices, we can formulate analogues of 
Propositions~\ref{propborderR}, \ref{propcritR}, \ref{theorealidealR}
for the radical ideal $I(\VC(I))$; the proofs of the first two results being similar are omitted.
  
\begin{proposition}\label{propborderC}
Let $t\ge d$ and $y\in \fCK_t$ 
for which $\rank \fCM_t(y)$ is maximum.
If, for some integer $1\le s\le t$,
\begin{equation}\label{flatiC}
\rank \fCM_s(y)=\rank M^\oC_{s-1}(y)=:r
\end{equation}
then $I(\VC(I))\supseteq \langle \Ker M^\oC_s(y)\rangle =:J$.
One 
can compute the set 
$W:=\VC(J)\supseteq \VC(I)$ which satisfies
 $\VC(I)\subseteq W$ and $|W|=r$.
Moreover, if $W=\VC(I)$ then $I(\VC(I))=J$.
\end{proposition}

\begin{proposition}\label{propcritC}
Let $t\ge d$ and $y\in \fCK_t$ for which $\rank \fCM_{t}(y)$ is maximum. 
Assume that, either
 (\ref{flatiC}) holds
 for some $2d\le s\le t$, or 
\begin{equation}\label{flatiiC}
\rank \fCM_{s}(y)=\rank M^\oC_{s-d}(y)
\end{equation}
for some $d\le s\le t$.
Then $I(\VC(I))=\langle \Ker M^\oC_s(y)\rangle$ (and one
 can find $\VC(I)$).
\end{proposition}

\begin{proposition} \label{theorealidealC}
Assume $|\VC(I)|<\infty$. 
\begin{description}
\item[(i)] If $V_\oC(I)=\emptyset$, then $ K^{2\oC}_{t}=\emptyset$ for $t$ large enough.
\item[(ii)] If $V_\oC(I)\ne \emptyset$ then, for $t$ large enough,
there exists $d\le s\le t$ such that 
$\rank \fCM_s(y)=\rank M^\oC_{s-d}(y)$ for all $y\in K^{2\oC}_{t}$.
\end{description}
\end{proposition}

\begin{proof}
Let $t\ge 2d$ and $y\in  K^{2\oC}_{t}$. Then, $h_1,\ldots,h_m\in \Ker M_t^\oC(y)$, since
$(M^\oC_t(y)h_j)_{\alpha'}=(h_jy)_{0\alpha'}=0$ for $|\alpha'|\le t\le 2t-2d_j$ by the assumption 
$\fCM_{t-d_j}(h_jy)=0$. Hence, $h_j\in \Ker \fCM_t(y)$ and thus 
$\overline{h_j}\in \Ker \fCM_t(y)$ too.

Let $\{g_1,\ldots,g_k\}$ be a Gr\"obner basis of $I(\VC(I))$ for a graded monomial ordering. 
Analogously to Claim \ref{claim1}, one can show, using Hilbert's Nullstellensatz, the existence of $t_0\in\oN$
for which $g_l\in \Ker \fCM_t(y)$ for all $l$ and $t\ge t_0$.
If $\VC(I)=\emptyset$, then $1\in \Ker \fCM_t(y)$ which implies $y_0=0$, thus showing 
$K^{2\oC}_{t}=\emptyset$. Assume now that $\VC(I)\ne \emptyset$.
Let $\BB$ be the basis of $\PC/I(\VC(I))$ for the chosen  monomial ordering
and set  $d_\BB:=\max_{b\in \BB} \deg(b)$  (which is well defined as $\BB\ne \emptyset$) and
$\CC:=\{\overline bb'\mid b,b'\in \BB\}$. 
Then any monomial $\ox^\alpha x^{\alpha'}$
can be written $\ox^\alpha x^{\alpha'}=r^{(\alpha\alpha')} +\sum_{l,l'=1}^k u^{(\alpha\alpha')}_{l l'}
\overline{g_l} g_{l'}$ where $u^{(\alpha\alpha')}_{l l'}\in \oC[x,\ox]$,
$r^{(\alpha \alpha')}\in \Span_\oC(\CC)$, and $\deg(u^{(\alpha\alpha')}_{l l'}\overline{g_l} g_{l'})
\le |\alpha|+|\alpha'|$.
Let $t_1:= \max(3d,d+2d_\BB)$ and $t\ge t_1+1$.
Then $\ox^\alpha x^{\alpha'} -r^{(\alpha \alpha')}\in \Ker \fCM_t(y)$
whenever $|\alpha+\alpha'|\le t_1$ which, together with
$\deg(r^{(\alpha \alpha')})\le 2d_\BB \le t_1-d$, shows that 
$\rank \fCM_{t_1}(y)=\rank \fCM_{t_1-d}(y)=:r$.

There remains to show that $\rank M^\oC_{t_1-d}(y)=r$.
Applying Theorems \ref{theoCFmeasure}, \ref{theoCFflat}, there exists $W\subseteq \oC^n$, $|W|=r$, $\lambda_v>0$ ($v\in W$)
such that, if we set  $\tilde y:=\sum_{v\in W}\lambda_v \zeta_{\ov}\otimes \zeta_v$,
then $\fCM(\tilde y)$ is a flat extension of $\fCM_{t_1}(y)$.
This implies $M^\oC_{t_1-d}(y)=\sum_{v\in W} \lambda_v  \zeta_{t_1-d,\ov}
\zeta_{t_1-d,v}^T$. 
As $h_j\in \Ker M^\oC_{t_1-d}(y)$ (since 
 $t_1 -d \ge \deg(h_j)$ as $t_1\ge 3d$),
we deduce that $h_j(v)=0$ for all $v\in W$ and thus $W\subseteq \VC(I)$. 
We now show that the vectors $\zeta_{t_1-d,v}$ ($v\in W$) are linearly independent, which implies that 
$\rank M^\oC_{t_1-d}(y) =|W|=r$, thus concluding the proof.
For this, consider interpolation polynomials $p_v\in \PC$ ($v\in W$), i.e., satisfying 
$p_v(w)=\delta_{v,w}$ for $v,w\in W$. One may assume that $\deg(p_v)\le d_\BB$
(replacing if necessary $p_v$ by its residue modulo $I(\VC(I))$ with respect to the basis $\BB$).
Assume $\sum_{v\in W}c_v \zeta_{t_1-d,v}=0$ for some $c_v\in \oC$; we show that all $c_v$'s are zero.
As $t_1-d\ge d_\BB$, we can take the scalar product with $\ve{p_w}$ ($w\in W$)
which yields $0=\sum_{v\in W}c_v p_w(v)=c_w$. 
$\hfill\Box$
\end{proof}

\medskip
\begin{remark}\label{comments}
Under condition  (\ref{flati}) or (\ref{flatiC}), the assumptions (i)-(iii) of 
Proposition \ref{propborder} hold.
Namely, one can construct an order ideal $\BB\subseteq \oT_{n,s-1}$ 
indexing a maximum nonsingular principal submatrix
of $M^\oK_{s-1}(y)$.
Moreover, one can choose  $\BB=\BB_\succ$, the set of standard monomials for the ideal
$J:= \langle \Ker M^\oK_s(y)\rangle$ with respect to a graded lexicographic order;
this is possible since $J$
 is zero-dimensional and there is a basis
in $\oT_{n,s-1}$ for $\PK/J$, which implies 
 $\BB_\succ\subseteq \oT_{n,s-1}$ by Lemma \ref{lemgrad}.
Such a basis $\BB_\succ$ can be found using the greedy sieve algorithm described in Section \ref{secsieve}. The execution of the 
algorithm requires checking 
whether some set $T\subseteq \oT_{n,s-1}$ is linearly independent in $\PK/J$. 
Using Corollary \ref{corker},
this can be checked by testing whether 
$T$ indexes a nonsingular 
principal submatrix of $M^\oK_t(y)$, thus by a rank computation on $M_t^\oK(y)$.
Finally the formal multiplication matrices as defined in Proposition \ref{propborder} coincide with the 
(usual) multiplication matrices in $\PK/J$ and thus they commute pairwise.

Therefore, the conclusion of Proposition \ref{propborder} also applies under (\ref{flati})
or (\ref{flatiC}):
If $W=V_\oK(I)$ then one can construct a border basis of $J=I(V_\oK(I))$.
Moreover, when using $\BB_\succ$, one can also construct the reduced Gr\"obner basis of $J$ 
for the graded lexicographic monomial ordering.
A sufficient condition for $W=V_\oK(I)$ is given above in Proposition~\ref{propcritR} ($\oK=\oR$) and Proposition~\ref{propcritC} ($\oK=\oC$). 

In fact, as will be explained in Section \ref{secgrobner}, we can adapt this strategy to find a Gr\"obner basis for an arbitrary monomial ordering.
\end{remark}

\begin{remark}\label{remrealideal}
The above results involve the matrix $\fCM_t(y)$ where the argument $y\in \oC^{\oN^{2n}}$ is a complex sequence 
satisfying (\ref{relher}). As explained in Section \ref{secmeasure} above,
one may work instead with the moment matrix $M^\oR(a)$ where $a\in \oR^{\oN^{2n}}$ is the real sequence 
defined as in (\ref{link}), and Propositions \ref{propborderC} and \ref{propcritC} could be reformulated in terms of real sequences only.

In fact, when  the ideal $I$ is generated by real polynomials $h_1,\ldots,h_m$,
 its set $\VC(I)$ 
of complex roots
is  closed under complex conjugations, i.e., $v\in \VC(I)$ $\Longleftrightarrow$ 
$\bar v\in \VC(I)$ and, for a polynomial $f\in \PC$, $f\in I(\VC(I))$ 
if and only if its real and imaginary parts belong to $I(\VC(I))$; that is,
 it suffices to determine $I(\VC(I))\cap\PR$.
For this, it suffices to consider {\em real valued} matrices 
 $M^\oC_t(y)$ (or $\fCM_t(y)$), i.e., with $y\in K^\oC_t\cap \oR^{\oN^{2n}_{2t}}$
(or $\fCK_t\cap \oR^{\oN^{2n}_{2t}}$)  in Propositions
  \ref{propborder}, \ref{propborderC} and \ref{propcritC}.
Indeed, one may e.g. easily verify that Lemma \ref{lemkerC} remains valid within the context of real polynomials and replacing $K^\oC_t$ (or $\fCK_t$) by $K^\oC_t\cap\oR^{\oN^{2n}_{2t}}$
(or $\fCK_t\cap \oR^{\oN^{2n}_{2t}}$). 
(Use here the fact that, since $\VC(I)$ is closed under conjugation, then
$\frac{1}{2}(\zeta_{2t,\bar v}\otimes \zeta_{2t,v} +\zeta_{2t,v}\otimes \zeta_{2t,\bar v})$
belongs to $\fCK_t\cap \oR^{\oN^{2n}_{2t}}$.)
\end{remark}

\subsection{Algorithm and implementation}
\label{subsec::algorithm}

With the results of Sections \ref{secfind1}-\ref{secfind3} 
 we have all the ingredients needed to compute the radical ideals $I(V_\oR(I))$ and
$I(\VC(I))$ of an ideal $I$ given by its generators. We now describe the algorithm in more detail.

 For convenience, let $K_t$
(resp., $M_t(y)$) stand for $K^\oR_t$, $K^\oC_t$,  $\fCK_t$ (resp.,
$M^\oR_t(y)$,  $M^\oC_t(y)$, $\fCM_t(y)$). For the task of computing $I(V_\oR(I))$,
we will use $K_t=K^\oR_t$ (and apply Propositions \ref{propborder}, 
\ref{propborderR}, \ref{propcritR})
and for the task of computing $I(\VC(I))$ we use $K_t=K^\oC_t$ (and apply Proposition \ref{propborder}) or $K_t=\fCK_t$
(and apply Propositions \ref{propborderC}, \ref{propcritC}).
The algorithm consists of five main parts:  For a given order $t\ge d$,
\begin{enumerate}
\item[{\bf (i)}] Find an element $y\in K_t$ maximizing the rank of $M_t(y)$.
\item[\bf (ii)] Check the ranks of the principal submatrices of $M_t(y)$.
\item[\bf (iii)] Compute a basis for the column space of $M_{s-1}(y)$ and the 
quotient space $\PK/J$ ($J=\langle \Ker M_{s}(y)\rangle$, for suitable 
$1\le s\le t$). 
\item[\bf (iv)] Compute  the formal multiplication matrices.
\item[\bf (v)] Construct a basis for the ideal $ J$.
\end{enumerate}
In step (ii) we search for a submatrix $M_s(y)$ of $M_t(y)$ 
satisfying Proposition~\ref{propborder} (i)-(iii), or the rank condition (\ref{flati}) (resp. (\ref{flatiC})), or (\ref{flatii}) (resp. (\ref{flatiiC})).
 Depending on what condition
is satisfied, 
the algorithm returns  a subideal $J\subseteq I(V_\oK(I))$ together with a superset $W\supseteq V_\oK(I)$, 
or the desired radical ideal $I(V_\oK(I))$ and the desired set of roots $V_\oK(I)$.
One can anyway verify a posteriori whether $W=V_\oK(I)$, simply by checking whether $h_j(v)=0$ for all $j\le m$ and  $v\in W$.
In the sequel of this section we give more details about these different tasks.

\subsubsection{Finding $y\in K_t$ maximizing the rank of $M_t(y)$}
\label{rem::SDP}
This first task
can be cast as the problem of finding a feasible solution of a semidefinite program, that has maximum rank. For details on the theory and applications of {semidefinite programming} the interested reader is referred, e.g., to \cite{VB95}, \cite{WoSaVa00}.
It is a known geometric property of semidefinite programs 
that a feasible solution 
has maximum rank if and only if it lies in the relative interior of the feasible region 
and that such point can be found with interior-point algorithms using self-dual embedding (see, e.g., \cite{K02}, \cite{WoSaVa00}).
Let us give some details.

\medskip
Consider a general instance of semidefinite program
\begin{equation}\label{P}
p^*:=\inf \sum_{j=1}^mb_jy_j\ \text{ s.t. } \sum_{j=1}^m y_jA_j-C \succeq 0
\end{equation}
and its dual semidefinite program
\begin{equation}\label{D}
d^*:=\sup \Tr(CX) \text{ s.t. } \Tr(A_jX)=b_j \ (j=1,\ldots,m), X\succeq 0.
\end{equation}
Here, $A_j,C,X$ are Hermitian matrices, $b,y\in \oR^m$, $X$, $y$ are the variables.
Obviously, $d^*\le p^*$ (weak duality). There is no duality gap (i.e., $p^*=d^*$),
 e.g., 
when (\ref{P}) is strictly feasible (i.e., $\exists y\in \oR^m$ with
$\sum_{j=1}^m y_jA_j-C \succ 0$) or when (\ref{D}) is strictly feasible 
(i.e., $\exists X\succ 0$ feasible for (\ref{D})).
When (\ref{D}) is strictly feasible and $d^*<\infty$, then (\ref{P}) attains its
minimum, i.e., the set of optimal solutions is nonempty.
The feasible region to (\ref{P}) is the convex set
$$ \begin{array}{ll}
K =\{y\mid \sum_{j=1}^m y_jA_j-C \succeq 0\}
 =\{y\mid 
u^*( \sum_{j=1}^m y_jA_j-C)u\ge 0 \ \forall u\in \oK^m\}.
\end{array}$$
Therefore, for $y\in K$, $y$  lies in the relative interior of $K$ 
if and only if $\Ker ( \sum_{j=1}^m y_jA_j-C)\subseteq \Ker ( \sum_{j=1}^m z_jA_j-C)$
for all $z\in K$ or, equivalently, if  $ \sum_{j=1}^m y_jA_j-C$ has maximum possible rank
(same argument as for Lemma~\ref{lemkerR}).

Semidefinite programs
can be solved in polynomial time to an arbitrary precision using, e.g.,
 the ellipsoid method, whose running time is however prohibitively high in practice.
 Interior-point methods are now the method of choice 
for solving semidefinite programs.
Assuming strict feasibility of (\ref{P}) and (\ref{D}), interior-point algorithms construct 
sequences of points on the so-called central path, which has the property of  
converging to an optimum solution of maximum rank \cite{GoSc98}. 
One can also find a maximum rank optimum solution under the weaker assumption that
(\ref{P}), (\ref{D}) are feasible (but not necessarily strictly feasible), if $p^*$ is attained,
     and $p^*=d^*<\infty$.
Indeed one can then construct the so-called 
extended self-dual embedding which is a strictly feasible semidefinite program with the property that a maximum rank optimum solution to it
yields a maximum rank optimum solution to the original problem (\ref{P}) (see e.g. \cite[Ch.~4]{K02}, \cite[Ch.~5]{WoSaVa00}).

\medskip
For our problem of finding $y\in K_t$ maximizing $\rank M_t(y)$,
consider the semidefinite program
\begin{equation}\label{P0}
p^*:=\min \ 1 \text{ s.t. } M_t(y)\succeq 0,\ M_{t-d_j}(h_jy)=0 \ (j=1,\ldots,m), y_0=1,
\end{equation}
where we add the condition (\ref{relher}) in the complex case. 
One can interpret (see, e.g., \cite{Las01}) the dual of (\ref{P0}) as 
\begin{equation}\label{D0}
\begin{array}{lc}
d^*:=\max \ \lambda  \text{ s.t. } & 1-\lambda = s +\sum_{j=1}^m q_j h_j 
 \text{ with } s,q_j \text{ polynomials }\\
&  \deg(s), \deg(q_jh_j)\le 2t, \
 s \text{ is s.o.s.}
\end{array}
\end{equation}
where `$s$ is s.o.s.' means that $s$ can be written as a sum of squares, i.e., 
$s=\sum_h |u_h|^2$ for some polynomials $u_h\in \PK$ or $\oC[x,\ox]$.
Obviously, (\ref{P0}) is feasible if $V_\oK(I)\ne \emptyset$. Moreover,
 (\ref{D0}) is feasible (e.g. with $\lambda=1$, $s=q_j=0$ as feasible solution) and,
if $K_t\ne \emptyset$, then
$p^*=1$ is attained by the whole set $K_t$ and
$p^*=d^*=1$. Hence 
an interior-point algorithm implementing the self-dual embedding technique applied to problem (\ref{P0}) is guaranteed to return the following information\footnote{Three options (I),(II),(III) are
described in \cite[Ch.~5, p. 119]{WoSaVa00}; (i) corresponds to (I) and (ii) to (II),(III). Indeed, under (III) a certificate is reported that
no complementary pair exists which implies, in our case, that (\ref{P0}) is infeasible, since
any $y\in K_t$ together with the solution $\lambda=1$, $s=q_j=0$ to (\ref{D0})
makes a complementary pair.}:
Either (i) $y\in K_t$ maximizing $\rank M_t(y)$, or (ii)
a certificate that (\ref{P0}) is infeasible thus implying $V_\oK(I)=\emptyset$.
For our computations we use the semidefinite programming solver SeDuMi-1.05 \cite{St99,St02}
which has this feature. Practically, this means that the solution returned by the algorithm is 
very close to a  maximum rank optimum solution. 
\begin{remark}
When using a semidefinite programming solver without the maximum rank property, one can recover a maximum rank solution 
to (\ref{P0})  from a feasible solution $\hat y$ to (\ref{P0}), using the following simple iterative algorithm.
Let $u_1,\ldots,u_p$ be a set of vectors that span $\Ker M_t(\hat y)$,
set $C:=\sum_{i=1}^pu_iu_i^*$, and consider the semidefinite program:
$\max \langle C,M_t(y)\rangle $ subject to $y$ satisfying the constraints of (\ref{P0}).
If the optimum value is equal to 0, then $\hat y$ is in fact a solution of maximum rank.
Otherwise, let $y_1$ be the optimum solution returned by the solver; then $\Ker M_t(\hat y)\not\subseteq \Ker M_t(y_1)$.
Then, $y_2:=\frac{1}{2}(\hat y+y_1)$ is feasible for (\ref{P0}) and
$\Ker M_t(y_2)=\Ker M_t(\hat y)\cap \Ker M_t(y_1) 
\subset \Ker M_t(\hat y)$. Hence we have found a feasible solution $y_2$ to (\ref{P0})
for which the rank of $M_t(y_2)$ is larger   than that $M_t(\hat y)$.
Iterate replacing $\hat y$ by $y_2$.
\end{remark}

\subsubsection{\bf Checking ranks of submatrices of $M_t(y)$}
Once a maximum rank matrix $M_t(y)$ is found, one has to  check if for some $1\le s\le t$ the conditions of Proposition \ref{propborder} (i)-(iii) hold, or if (\ref{flati}) (resp. (\ref{flatiC})) holds, or if (\ref{flatii}) (resp. (\ref{flatiiC})) holds.
For this one has to compute the ranks of the  principal submatrices 
$M_s(y)$ of $M_t(y)$ for $s\leq t$. Checking the rank of a matrix consisting of numerical values is computationally sensitive. This is carried out using  singular value decomposition which at the same time can be used to generate a basis of the column space; see the next section for more details. The determination of the rank is done by detecting zero singular values or a decay of more than $1\text{e-3}$ between two subsequent values,
where  singular values less than $1\text{e-8}$ are declared to be zero.

\subsubsection{\bf Computing a basis for the column space of $M_{s-1}(y)$ and the
 quotient space $\PK/J $}

We indicate here how to compute a basis of the column space of the matrix $M_{s-1}(y)$.
Under some conditions (recall Corollary \ref{corker}), such basis also yields a basis of 
the quotient space  $\oK[x]/J$ (as before, $J:=\langle \Ker M_s(y)\rangle$), 
which is needed for the computation of the multiplication matrices.
The choice of this basis will have an influence on the numerical stability of the extracted 
set $W$ of  solutions and on the properties of the basis for $J$ as well.

\medskip\noindent
{\bf \textit{Using singular value decomposition.}}
It is a well known fact from linear algebra that a numerically stable way of finding an orthonormal basis $\BB$ for the column space of a matrix $M$ is to use its singular value decomposition (SVD):
$M = U\Sigma V^*,$
where $U,V$ are unitary and $\Sigma$ is diagonal with nonnegative entries.
The diagonal entries of $\Sigma$ are the singular values of $M$ (i.e., the square roots 
of the eigenvalues of $MM^*$); the number of nonzero diagonal entries of $\Sigma$
is thus equal to $r:=\rank \ M$.  Then the set $\{U_{1},\ldots,U_{r}\}$ of columns of $U$ corresponding to the nonzero diagonal entries of $\Sigma$ 
forms an orthonormal  basis of the column space of $M$. 
As we already did perform a SVD to determine the rank of the matrix $M:=M_{s-1}(y)$, this computation comes with no extra effort. 
For $i=1,\ldots,r$, let  $b_i:= \zeta_{s-1,x}^TU_{i}$ be the polynomial
with vector of coefficients $\ve{b_i}=U_{i}$. The next lemma shows that, under some rank condition,  
$\{b_1,\ldots,b_r\}$ is a basis of $\PK/J$. 

\begin{lemma}\label{lemSVD}
Let $\{U_1,\ldots,U_r\}$ be a basis of the column space of $M_{s-1}(y)$,
let $b_i:=\zeta_{s-1,x}^TU_i$ ($i=1,\ldots,r$), 
and assume 
that the rank condition (\ref{flati}) or (\ref{flatiC}) holds. Then
the set $\{b_1,\ldots,b_r\}$ is a basis of the quotient space
$\PK/J$.
\end{lemma}

\begin{proof}
As the rank condition (\ref{flati}) or (\ref{flatiC}) holds, we know from Corollary~\ref{corker}
that
$\dim \PK/J=r:=\rank M_{s-1}(y)$ and 
 $J\cap \oK_{s-1}[x]= \Ker M_{s-1}(y)$.
Hence it suffices to show that $\{b_1,\ldots,b_r\}$ is linearly independent in 
 $\PK/J$.
For this assume $\sum_{i=1}^r\lambda_i b_i\in J$, i.e., $\zeta_{s-1,x}^T(\sum_{i=1}^r\lambda_i U_i) \in J$.
The vector $p:=\sum_{i=1}^r\lambda_i U_i $ lies in the column space of $M_{s-1}(y)$.
On the other hand, $p\in \Ker M_{s-1}(y)$ since the corresponding polynomial 
$p$ lies in $J\cap \oK_{s-1}[x]$. Therefore, $p=0$ which implies that all $\lambda_i=0$.
$\hfill\Box$ \end{proof}

\medskip\noindent
{\bf \textit{Using a `greedy' algorithm.}}
If we want to compute a border basis or a Gr\"obner basis with our algorithm, we need a monomial basis $\BB$ (i.e.,
$\BB\subseteq \oT_n$)  for the quotient space $\oK[x]/J$.
For this, it suffices to construct a set $\BB$
indexing a maximum principal nonsingular submatrix of $M_{s-1}(y)$ as
$\BB$ 
is then a basis of $\PK/J$ under the conditions of Corollary~\ref{corker}.
One can apply the following simple procedure (proposed in \cite{Lau2})
 for constructing $\BB$:
Scan monomials in $\oT_{n,s-1}$ by increasing degree, starting with $t_0=1,t_1=x_1,t_2=x_2, \ldots$. 
Initialize $\BB:=\{t_0\}$.
Let  $\BB$ be the current set and $t_k$ be the current monomial to be scanned.
If $\BB\cup \{t_k\}$ indexes a linearly independent set of columns of $M_{s-1}(y)$, then reset $\BB:=\BB\cup\{t_k\}$, otherwise scan the next monomial
$t_{k+1}$. This procedure is `greedy' in the sense that one keeps
 adding as many low degree monomials as possible to the basis.
One can stop as soon as $|\BB|=r$.
Alternatively, one may construct a reduced row echelon form of $M_{s-1}(y)$ using Gauss Jordan elimination with partial pivoting (pivot variables serve as basis $\BB$), see~\cite{HeLa05}.
One can verify afterwards whether the constructed basis $\BB$ is an order ideal; it turns out that this is the case in most tested instances.

The greedy sieve algorithm  described earlier in Section \ref{secsieve}  
produces directly an order ideal basis.
Indeed, given any graded monomial ordering $\succ$, we  can apply it
to obtain the set $\BB=\BB_\succ$ of standard monomials for this ordering,
forming an order ideal basis of $\PK/J$  (as we know from Lemma \ref{lemgrad} that
$\BB_\succ$ is contained in $\oT_{n,s-1}$ under the conditions of Corollary \ref{corker}.)
See Section~\ref{secgrobner} for an extension to the case of an arbitrary monomial ordering.

Note that, although desirable from an algebraic point of view, monomial bases for $\oK[x]/J$ sometimes lead to a less accurate set $W$ of extracted solutions as compared to those extracted with a polynomial basis $\BB$ based on SVD;
 see e.g. Examples~\ref{katsura5}, \ref{philipp}.

\subsubsection{Computing formal multiplication matrices}
Let $\BB=\{b_1,\ldots,b_r\}$ be a basis of the column space of $M_{s-1}(y)$.
By the assumptions of 
Proposition \ref{propborder} or under the rank conditions (\ref{flati}) or (\ref{flatiC}),
there exist scalars $a^{(ij)}_k$ ($k=1,\ldots,r$) for which 
$x_ib_j-\sum_{k=1}^r a^{(ij)}_kb_k\in \Ker M_s(y)$, for all $i=1,\ldots,n$, $j=1,\ldots,r$.
Then the vector $(a^{(ij)}_k)_{k=1}^r$ is the $j$th column of the (formal)
multiplication matrix $\XX_{x_i}$.
  We indicate how to compute $\XX_{x_i}$ from $M_s(y)$.

Suppose first that $\BB$ is a monomial basis, i.e., $\BB\subseteq \oT_{n,s-1}$.
Let $M_\BB$ denote the principal submatrix of $M_s(y)$ indexed by $\BB$ and let $P_{x_i}$ be 
the submatrix of $M_s(y)$ whose 
rows are indexed by $\BB$ and whose columns are indexed by the set
$x_i\BB:=\{x_ib_j\mid j=1,\ldots,r\}$.
As observed in \cite{Lau2}, we have 

\begin{equation}\label{multmono}
\XX_{x_i}= M_\BB^{-1}P_{x_i}.
\end{equation}
Indeed, for $b\in \oT_{n,s}$, let $C_b$ denote the column of $M_s(y)$ indexed by $b$ restricted to the rows indexed by $\BB$.
Then, 
\begin{equation}
\label{eqn::mulmatrix_mon}
C_{x_ib_j}=\sum_{k=1}^r a^{(ij)}_k C_{b_k}= M_\BB a^{(ij)}\,,
\end{equation} i.e.,
$a^{(ij)}=M_\BB^{-1} C_{x_ib_j}$, which gives
$\XX_{x_i}=M_\BB^{-1} P_{x_i}$.

\medskip
Suppose now that $\BB$ is a polynomial basis obtained via SVD, as explained above. That is,
$b_i=\zeta_{s-1,x}^TU_{i}$ where $\{U_{1},\ldots,U_{r}\}$ is an orthonormal basis 
of the column space of $M_{s-1}(y)$ and thus of $M_s(y)$ under the rank condition 
(\ref{flati}) or (\ref{flatiC}).
As in the monomial case, the formal multiplication matrices can be derived from 
$M_s(y)$.
Let $\tilde P_{x_i}$ denote the submatrix of $M_s(y)$ with columns indexed by $x_i \oT_{n,s-1}$ and 
with rows indexed by $\oT_{n,s-1}$. 
Let $U$ denote the matrix with columns $U_1,\ldots,U_r$, and set
$P_{x_i}:= U^T \tilde P_{x_i} U$ and 
$M_\BB:= U^T M_{s-1}(y) U$. Then, $M_\BB$ is nonsingular.
Moreover,
\begin{equation}\label{multSVD}
M_\BB \XX_{x_i}= P_{x_i},
\end{equation}
which allows the computation of $\XX_{x_i}=M_\BB^{-1}P_{x_i}=\Sigma^{-1} P_{x_i}$ where $\Sigma$ is the diagonal matrix containing the positive singular values of $M_{s-1}(y)$.
We verify that (\ref{multSVD}) holds.
By construction, the polynomial $x_ib_j-\sum_{k=1}^r a^{(ij)}_kb_k=
x_i\zeta_{s-1,x}^TU_j-\sum_{k=1}^r a^{(ij)}_k\zeta_{s-1,x}^TU_k$ lies 
in $\Ker M_{s}(y)$. This implies $0=\tilde P_{x_i}U_j-M_{s-1}(y)Ua^{(ij)}$
and thus $U^T\tilde P_{x_i}U_j= U^TM_{s-1}(y)Ua^{(ij)}= M_\BB a^{(ij)},$
which shows that the two matrices $P_{x_i}$ and 
$M_\BB\XX_{x_i}$ have identical $j$th columns.

\subsubsection{\bf Constructing a basis for the ideal $J:= \langle \Ker M_s(y)\rangle$}\label{secgrobner}

\medskip\noindent
{\bf \textit{A linear  basis of $\Ker M_s(y)$.}}
The simplest way of producing  a basis for the ideal $J= \langle \Ker M_s(y)\rangle$ is 
simply by considering a linear basis of $\Ker M_s(y)$. Such a basis can be found by 
using again a singular value decomposition for $M_s(y)$. Indeed, if $M_s(y)=U\Sigma V^*$ is the SVD, 
then the columns $V_i$  of $V$ corresponding to the zero diagonal entries of $\Sigma$ (the zero singular values of $M_s(y)$) form an orthonormal basis of
$\Ker M_s(y)$.
Then the polynomials $\zeta_{s,x}^TV_{i}$ corresponding to the zero singular values
of $M_s(y)$ form a basis of $J$.
A drawback of this basis however is that it is usually highly overdetermined and has a large cardinality, equal to $|\oT_{n,s}|-\rank M_s(y)$.

\medskip\noindent
{\bf \textit{A border basis.}}
As shown in~\cite[Sec. 8.2, Ch. 10]{St04}, 
it is desirable to avoid overdetermined bases for $J$ because
it could lead to inconsistencies in the basis for numerical reasons.
To avoid this drawback, border bases are proposed in \cite{St04}
and their numerical properties are investigated.
If during the construction of the formal multiplication matrices an order ideal  
basis $\BB$ of $\oK[x]/J$ was used, we deduce immediately a border basis 
consisting of the polynomials  
\begin{equation}
\label{eqn::constrBB}
x_ib_j-\sum_{k=1}^r a^{(ij)}_kb_k \ \text{ for }x_ib_j\in \pBB.
\end{equation}

\medskip\noindent
{\bf \textit{A  Gr\"obner basis.}}
If the monomial basis $\BB$ of $\PK /J$ is the set of standard monomials $\BB_\succ$ 
with respect to
a monomial ordering $\succ$ obtained, e.g., with the greedy sieve algorithm,
then  the border basis in (\ref{eqn::constrBB}) is actually a Gr\"obner basis with respect to
the monomial ordering $\succ$.
When $\succ$ is a graded monomial ordering then, in view of Lemma \ref{lemgrad},
$\BB_\succ\subseteq \oT_{n,s}$ and thus $\BB_\succ$ can be found with Algorithm 1 applied to 
$(I(V_\oK(I)),\succeq,s)$, using the following independence oracle: a subset of $\oT_{n,s}$ is 
independent in
$\PK/I(V_\oK(I))$ if and only if it indexes an independent set of columns of $M_s(y)$.
In general when $\succ$ is not a graded monomial ordering we are not assured to find $\BB_\succ$ within 
$\oT_{n,s}$. However we can proceed as follows. As $M_s(y)$ is a flat extension of $M_{s-1}(y)$, 
by the results in Section \ref{secflat}, there exists an extension $\tilde y\in\oR^{\oN^n_t}$ (for any $t\ge s$)
such that $M_t(\tilde y)$ is a flat extension of $M_s(y)$. 
As $\langle \Ker M_s(y)\rangle=
\langle \Ker M_{t}(\tilde y)\rangle = I(V_\oR(I))$, a subset of $\oT_{n,t}$ is independent in 
$\PK/I(V_\oK(I))$ if and only if it indexes an independent set of columns of $M_t(\tilde y)$.
Thus to find $\BB_\succ$, apply Algorithm 1 iteratively to $t=s+1,s+2,\ldots$ until finding 
$\BB_t=\BB_{t+1}$, in which case we know from Lemma \ref{lemalg} (ii) that $\BB_\succ=\BB_t$. 
Remains only to address the question on how to find the flat extension $\tilde y$. The existence proof in \cite{CF96} 
is constructive and can roughly be sketched as follows (see also \cite{Lau2} for details).
Say we want to construct a flat extension $C:=M_{s+1}(\tilde y)$ of $B:= M_s(y)$, under the assumption that
$B$ is a flat extension of $M_{s-1}(y)$. We indicate how to construct the column 
$C(\cdot,\gamma)$ of $C$ indexed by a monomial
$x^\gamma$ of degree $s+1$. Write, say, $x^\gamma=x_ix^\beta$. By assumption, the column 
$B(\cdot,\beta)$ of $B$ indexed by $x^\beta$ can be expressed as a linear combination
$\sum_{|\alpha|\le s-1}\lambda_\alpha B(\cdot,\alpha)$ of columns indexed by $\oT_{n,s-1}$;
then define $C(\cdot,\gamma)$ as $\sum_{|\alpha|\le s-1}\lambda_\alpha C(\cdot, \alpha+e_i)$.
Note that this construction relies on the fact that the kernel of moment matrices enjoys ideal-like properties.

\medskip\noindent
\subsubsection{\bf Summary of the algorithm}

Algorithm \ref{alg::borderbasis}
below summarizes our algorithm.
This algorithm has been implemented in Matlab using the Yalmip toolbox~\cite{Lo04}.
For solving the semidefinite program (\ref{P0}) the semidefinite solver 
SeDuMi-1.05 \cite{St99,St02} is used. 
As described above and can be seen in the examples in the next section, the rank detection is the most critical task.
This was the main motivation for the  weaker conditions from Section \ref{secfind1},
which extend the possibility of extracting solutions. In the examples below 
this deficiency is clearly indicated in some rank sequences not exactly matching the theory.

\begin{algorithm}
\caption{\emph{Numerical border basis computation:}}
\begin{algorithmic}[1]
\Require Polynomial generators $h_i$ for $I:=\langle
h_1,\ldots h_m\rangle \subseteq \oK[x]$ and relaxation order $t\in \oN$
\Ensure A basis for an ideal $J \subseteq I(V_\oK(I))$,  the set $\VC(J)$,
and a basis $\BB$ for the quotient ring $\oK[x] / J$ 
\State Solve the SDP (\ref{P0}). If the SDP is infeasible, return $V_\oK(I)=\emptyset$. 
Otherwise, return a feasible solution $y$ for which $M_t(y)$ has maximum rank
\State Compute SVD for all principal submatrices $M_s(y)$ ($s=1,\ldots, t$)
\State Determine $\rank M_s(y)$ ($s=1,\ldots,t$)  and check whether 
the conditions of Prop.~\ref{propborder} -- \ref{propcritC} hold
\State Fix $s$ (for which one of Prop.~\ref{propborder} - \ref{propcritC} applies) 
\State Compute a basis $\BB$ of the column space of $M_{s-1}(y)$:
\Statex a) using the SVD decomposition ($\BB$ is a polynomial basis)
\Statex b) using a greedy algorithm ($\BB$ is a monomial basis)
\Statex c) using a greedy sieve algorithm ($\BB$ is the set of standard monomials for a monomial ordering $\succ$)
\State Compute the multiplication matrices $\XX_{x_i}=M_\BB^{-1} P_{x_i}$
\State Compute a basis for the ideal $J$
\Statex a) a SVD basis of $\Ker M_s(y)$ (requires $\rank M_s(y)=\rank M_{s-1}(y)$)
\Statex b) a border basis of $J$ (requires that $\BB$ is a monomial basis)
\Statex c) a Gr\"obner basis (requires that $\BB$ is the corresponding set of standard monomials)
\If{ the conditions of Prop.~\ref{propborder} are met}
\State \Return a border/Gr\"obner basis of $J \subseteq I(V_\oK(I))$,
 a basis $\BB$ of  $\oK[x] / J$, and the set $\VC(J)$
\Else
\State \Return ERROR: "No extraction possible. Increase relaxation order $t$."
\EndIf
\end{algorithmic}
\label{alg::borderbasis}

\begin{remark}
In view of Propositions~\ref{theorealidealR}, \ref{theorealidealC}, the algorithm terminates for $t$
large enough and finds $J=I(V_\oR(I))$ or $I(\VC(I))$.
The algorithm can also be used for testing existence of solutions. Let us give some details e.g. in the real case.
If at step 1 one detects infeasibility of the SDP then one can already conclude 
$V_\oR(I)=\emptyset$. Suppose now the SDP is feasible.
At step 4, as observed in Remark \ref{remdetect},
the conditions of Proposition \ref{propcritR} cannot be met if $V_\oR(I)=\emptyset$, but it could be that the conditions 
of Proposition \ref{propborder} or \ref{propborderR} are met. In that case one can extract a set
$W\supseteq V_\oR(I)$. Then one can simply test whether the points of $W$ satisfy the 
given equations $h_1=\ldots=h_m=0$ to detect whether $V_\oR(I)$ is empty or not.

\end{remark}
\end{algorithm}

\section{Numerical Examples}
\label{sec::numexamples}

We present here the results of our algorithm applied to some examples, mainly taken from the literature.
In each example, we specify the ideal $I$ by its generators $h_1,\ldots,h_m$.
\noindent
Let us explain Tables 1-5 below. At a given order $t$, let $y$ be the 
optimal solution to (\ref{P0}) returned by the SDP 
solver.
The abbreviations `MON' and `SVD' refer to using a {\em monomial} base 
of the quotient space, or a base found via the {\em SVD method}. \\
$\bullet$ The column `rank sequence' shows $(\rank M_0(y), \ldots,\rank M_t(y))$. \\
$\bullet$ The column 
`extract. order' shows some numbers 
$s_{\text{mon}}(r_{\text{mon}})/s_{\text{svd}}(r_{\text{svd}})$.
When using a monomial base, 
$r_{\text{mon}}$ is the smallest order at which the extraction procedure could be carried out and
$s_{\text{mon}}$ is the order at which it was effectively carried out and gave the results reported here; analogously with the SVD method.\\
$\bullet$ The column `accuracy' shows the accuracy of the returned solutions, 
i.e., $\max_{j,x} |h_j(x)|$, where $h_j$ runs over the generators of $I$ and $x$ over the extracted solutions.\\
$\bullet$ The column `comm. error' shows the commutativity error for the multiplication matrices, i.e.,
$\max_{i,j=1}^n \text{abs}(\XX_{x_i}\XX_{x_j}-\XX_{x_j}\XX_{x_i})$ (where $\text{abs}(M)$ is the maximum absolute value of the entries of a matrix $M$). If the parameter `comm. error' is more than 1e-2, the multiplication matrices do not commute sufficiently and we then do not extract
solutions.

\begin{example}\label{cox3}
This simple example  from \cite[p.40]{CLO98} has two roots, both real.
\begin{align*} h_1&= x_2^4 x_1+3 x_1^3-x_2^4-3 x_1^2\\
h_2&=x_1^2 x_2-2 x_1^2\\
h_3&=2 x_2^4 x_1-x_1^3-2 x_2^4+x_1^2
\end{align*}
\begin{table}[ht!]
\begin{tabular}{|c|c|c|c|c|}
\hline 
order  & rank sequence & extract. order &  accuracy  & comm. error \\
$t$          &               &      MON/SVD     &  MON/SVD   & MON/SVD  \\
\hline 
\hline 
3&1  3  5  9&    ---    &   ---   & ---  \\
\hline 
4&1  2  2  2  7&4(2)/3(2)&1.9717e-9/0.00013144&9.676e-10/3.3908e-6\\
\hline 
5&1  2  2  2  2  8&4(2)/4(2)&2.9557e-8/3.5325e-5&1.8781e-11/1.2291e-6\\
\hline 
\end{tabular}
\caption{Results for Example~\ref{cox3}}
\end{table}

\noindent
Monomial basis of $\PR/I(V_\oR(I))$:
\begin{align*}\BB=\{1,x_1\}.\end{align*}
Border basis for $I(V_\oR(I))$ (showing in {\bf bold} the monomials in $\partial \BB$):
\begin{align*}
g_1&=-x_1+{\bf x_1^2},\nonumber\\
g_2&=-2 x_1+{\bf x_2},\nonumber\\
g_3&=-2 x_1+{\bf x_1 x_2}.
\end{align*}
Extracted real solutions $V_\oR(I)$:
\begin{align*}
x_1&=(2.12\text{e-8},1.91\text{e-6}),\\
x_2&=(1,2).
\end{align*}
The first two polynomials $g_1,g_2$ of the extracted border basis form a 
reduced Gr\"obner basis with respect to the graded reversed term order with 
 $x_1\prec x_2$.
The basis of $\sqrt I$ ($=I(V_\oR(I))$ as all roots are real) given in 
\cite{CLO98} has the form:
\begin{align*}
\{ &x_2^4  x_1 + 3 x_1^3  - x_2^4  - 3 x_1^2 , x_1^2  x_2 - 2 x_1^2 ,\nonumber\\ &2 x_2^4  x_1 - x_1^3  - 2 x_2^4  + x_1^2 , x_1 (x_1 - 1), x_2 (-2 + x_2)\}
\end{align*}
and is obtained via Seidenberg's method described in the paragraph `Related literature' in the Introduction.
Computing a graded reversed Gr\"obner Basis of $\sqrt I$ (using \texttt{tdeg} in \textsc{Maple})
leads again to the set $\{g_1,g_2\}$ found by our method. Thus our method finds 
here a simpler set of generators for $\sqrt I$ than the classical method of Seidenberg.
\end{example}

\begin{example}
\label{bifurcation}
This example is taken from the polynomial testsuite \cite{BM96}
(see \verb1http://www-sop.inria.fr/saga/POL/BASE/2.multipol/bifurc.html1). 
It has 20 complex solutions among which 8 are real. This example
illustrates the possibility of extracting solutions based on 
Proposition~\ref{propborder} in case none of the rank conditions are satisfied. 
\begin{align*} h_1&=5 x_1^9-6 x_1^5 x_2+x_1 x_2^4+2 x_1 x_3 \\
h_2&= -2 x_1^6 x_2+2 x_1^2 x_2^3+2 x_2 x_3 \\
h_3&= x_1^2+x_2^2-0.265625 
\end{align*}

\begin{table}[ht!]
\begin{tabular}{|c|c|c|c|c|}
\hline 
order  & rank sequence & extract. order &  accuracy  & comm. error \\
$t$          &               &      MON/SVD     &  MON/SVD   & MON/SVD  \\
\hline 
\hline 
5&1   4   8  16  25  34&    ---    &   ---   & ---  \\
\hline 
6&1   3   9  15  22  26  32&    ---    &   ---   & ---  \\
\hline 
7&1   3   8  10  12  16  20  24&3(3)/---(---)&0.12786/---&0.00019754/---\\
\hline 
8&1   4   8   8   8  12  16  20  24&4(3)/3(3)&4.6789e-5/0.00013406&4.7073e-5/0.00075005\\
\hline 
\end{tabular}
\caption{Results for Example~\ref{bifurcation}}
\end{table}

\noindent
Monomial basis of $\oR[x]/I(V_\oR(I))$:
\begin{align*}\BB=\{1,x_1,x_2,x_3,x_1^2,x_1 x_2,x_1 x_3,x_2 x_3\}\,.\end{align*}
Border basis of $I(V_\oR(I))$: 
\begin{align*}
g_1&=-0.28479 x_1+0.44124 x_1 x_2-1.5403 x_1 x_3+ {\bf x_1^3}\,,\\
g_2&=-1.7276 x_3-0.080949 x_1^2+8.1433 x_2 x_3+{\bf x_1^2 x_2}\,,\\
g_3&=-0.28763 x_3-0.0010314 x_1^2+0.48126 x_2 x_3+{\bf x_1^2 x_3}\,,\\
g_4&=-0.0015073 x_1+0.01299 x_1 x_2-0.12111 x_1 x_3+{\bf x_1 x_2 x_3}\,,\\
g_5&=-0.26563+x_1^2+{\bf x_2^2}\,,\\
g_6&=0.019164 x_1-0.44124 x_1 x_2+1.5403 x_1 x_3+{\bf x_1 x_2^2}\,,\\
g_7&=0.022008 x_3+0.0010314 x_1^2-0.48126 x_2 x_3+{\bf x_2^2 x_3}\,,\\
g_8&=-0.0018637 x_3-0.00067043 x_1^2+0.026066 x_2 x_3+{\bf x_3^2}\,,\\
g_9&=-0.00015166 x_1+0.00025958 x_1 x_3+{\bf x_1 x_3^2}\,,\\
g_{10}&=-0.0017335 x_3+0.01615 x_2 x_3+{\bf x_2 x_3^2}\,.
\end{align*}
Extracted real solutions:
\begin{align*}
x_1&=(-0.515,-0.000153,-0.0124)\,,\\
x_2&=(-0.502,0.119,0.0124)\,,\\
x_3&=(0.502,0.119,0.0124)\,,\\
x_4&=(0.515,-0.000185,-0.0125)\,,\\
x_5&=(0.262,0.444,-0.0132)\,,\\
x_6&=(-2.07\text{e-5},0.515,-1.27\text{e-6})\,,\\
x_7&=(-0.262,0.444,-0.0132)\,,\\
x_8&=(-1.05\text{e-5},-0.515,-7.56\text{e-7})\,.
\end{align*}
\end{example}

\begin{example} \label{gauss}
We now give an example for finding  $I(V_\oR(I)\cap S)$, where
 $I=\langle h_1,\ldots,h_4\rangle$ with
 \begin{align*} h_1& = x_1 + x_2 - 2 \,,\\
 h_2&=x_1x_3 + x_2 x_4 \,,\\
 h_3&=x_1x_3^2 + x_2 x_4^2 - \frac{2}{3}  \,,\\
 h_4&=x_1 x_3^3 + x_2 x_4^3 
\end{align*}
and $S$ is defined by the polynomial inequalities
$-1 \le x_1,x_2,x_3,x_4\le 1$.
This example, taken from \cite{VeGa95}, 
represents a Gaussian quadrature formula with 2 
weights and 2 knots, where   one is interested only in 
the  roots lying in the box 
$[-1,+1]$. 

\begin{table}[ht!]
\begin{tabular}{|c|c|c|c|c|}
\hline 
order $t$ & rank sequence & extract. order &  accuracy  & comm. error \\
          &               &      MON/SVD     &  MON/SVD   & MON/SVD  \\
\hline 
\hline 
2&1   2  11&2(2)/---(---)&0.00010224/---&1.1124e-9/---\\
\hline 
3&1   2   2  18&2(2)/2(2)&1.8985e-14/5.1015e-14&1.2212e-15/1.4155e-15\\
\hline 
4&1   2   2   2  24&2(2)/2(2)&3.5527e-15/8.5487e-15&2.2204e-16/1.1102e-16\\
\hline 
\end{tabular}
\caption{Results for Example~\ref{gauss}}
\label{tab::ex_gauss}
\end{table}

\noindent
Monomial basis of $\PR/I(V_\oR(I)\cap S)$:
\begin{align*}\BB=\{1,x_3\}.\end{align*}
Border basis of $I(V_\oR(I)\cap S)$:

\begin{align*}
g_1&=-1+{\bf x_1}\\
g_2&=-x_3+{\bf x_1 x_3}\\
g_3&=-1+{\bf x_2}\\
g_4&=-x_3+{\bf x_2 x_3}\\
g_5&=-0.33333+{\bf x_3^2}\\
g_6&=x_3+{\bf x_4}\\
g_7&=0.33333+{\bf x_3 x_4}
\end{align*}
For this example the border basis  is in fact a Gr\"obner basis, e.g. with respect to a graded lexicographic order with $x_1 \succ x_2 \succ x_4 \succ x_3$. Note however that it is not a reduced Gr\"obner basis since $g_2,g_4$ and $g_7$ are redundant for this ordering.

\noindent
Extracted real solutions $V_\oR(I)\cap S$:
\begin{align*}
x_1&=(1,1,-0.577,0.577)\,,\\
x_2&=(1,1,0.577,-0.577)\,.
\end{align*}
\end{example}

\begin{example}
\label{katsura5}
This example: Katsura 5, is an example in $\oR^6$
with 32 complex roots, including
12 real roots.
It is     taken from \\
\verb1http://www.mat.univie.ac.at/~neum/glopt/1
\verb1coconut/Benchmark/Library3/katsura5.mod1.
\begin{align*} h_1&=2 x_6^2+2 x_5^2+2 x_4^2+2 x_3^2+2 x_2^2+x_1^2-x_1\,,\\
h_2&=x_6 x_5+x_5 x_4+2 x_4 x_3+2 x_3 x_2+2 x_2 x_1-x_2 \,,\\
h_3&=2 x_6 x_4+2 x_5 x_3+2 x_4 x_2+x_2^2+2 x_3 x_1-x_3\,,\\
h_4&=2 x_6 x_3+2 x_5 x_2+2 x_3 x_2+2 x_4 x_1-x_4 \,,\\
h_5&=x_3^2+2 x_6 x_1+2 x_5 x_1+2 x_4 x_1-x_5 \,,\\
h_6&=2 x_6+2 x_5+2 x_4+2 x_3+2 x_2+x_1-1\,.
\end{align*}

\begin{table}[ht!]
\begin{tabular}{|c|c|c|c|c|}
\hline 
order  & rank sequence & extract. order &  accuracy  & comm. error \\
$t$          &               &      MON/SVD     &  MON/SVD   & MON/SVD  \\
\hline 
\hline 
1&1  7&    ---    &   ---   & ---  \\
\hline 
2&1   6  16&    ---    &   ---   & ---  \\
\hline 
3&1   6  12  12&---/3(3)&---/1.1928e-005&---/2.3073e-007\\
\hline 
\end{tabular}\caption{Results for Example~\ref{katsura5}}
\label{tab::ex_katasura5}
\end{table}

\noindent
We cannot extract solutions with a monomial base since the multiplication matrices do not commute, but we can extract the following
real solutions using a SVD basis:
\begin{align*}
x_1&=(0.277,0.226,0.162,0.0858,0.0115,-0.124)\,,\\
x_2&=(0.59,0.0422,0.327,-0.0642,-0.0874,-0.0132)\,,\\
x_3&=(1,-2.8\text{e-7},4.7\text{e-7},8.81\text{e-7},-2.79\text{e-6},-3.69\text{e-6})\,,\\
x_4&=(0.239,0.0608,-0.0622,-0.0233,0.186,0.219)\,,\\
x_5&=(0.441,0.151,0.0225,0.219,0.0935,-0.207)\,,\\
x_6&=(0.726,-0.0503,0.122,0.164,0.11,-0.208)\,,\\
x_7&=(0.462,0.309,0.0553,-0.102,-0.0844,0.0917)\,,\\
x_8&=(0.292,-0.101,0.181,-0.0591,0.193,0.141)\,,\\
x_9&=(0.753,0.0532,0.191,-0.114,-0.146,0.139)\,,\\
x_{10}&=(0.409,-0.0732,0.0657,-0.127,0.252,0.178)\,,\\
x_{11}&=(0.68,0.266,-0.154,0.0323,0.0897,-0.0735)\,,\\
x_{12}&=(0.136,0.0428,0.0417,0.0404,0.0964,0.211)\,,
\end{align*}
\end{example}

\begin{example}
\label{philipp}
 The following example  shows the  
 computation of the radical ideal using  complex moment matrices.
 \begin{align*} h_1& = x_1^2+x_2+x_3+1 \,,\\
 h_2 &=x_1+x_2^2+x_3+1\,,\\
 h_3&=x_1+x_2+x_3^2+1\,.
\end{align*}
This ideal is not radical and admits 7 solutions, among which the solution $(-1,-1,-1)$ has multiplicity two.
We solve the SDP program based on the set $\fCK_t$ (thus using {\em full} complex moment matrices). The rank sequence for full and pruned matrices $M^\oC_s(y)$ and $M^{2\oC}_s(y)$ are shown in Table~\ref{tab::ex_philipp}.

\begin{table}[ht!]
\begin{tabular}{|c|c|c|c|c|}
\hline 
order & rank sequence & extract. order &  accuracy  & comm. error \\
$t$& $M^\oC_s(y)$, ($M^{2\oC}_s(y)$)  &      MON/SVD     &  MON/SVD   & MON/SVD  \\
\hline 
\hline 
1&(1  4),&    ---    &   ---   & ---  \\
 &(1  7)&           &         &      \\
\hline 
2&(1  4  7), &    ---    &   ---   & ---  \\
&(1  7  7)&           &         &      \\
\hline 
3&(1  4  7  7), &3(3)/3(3)&0.0005719/0.00022538&0.00041241/0.00043871\\
& (1  7  7  7)&           &         &      \\
\hline 
\end{tabular}
\caption{Results for Example~\ref{philipp}}
\label{tab::ex_philipp}
\end{table}

\noindent
Monomial basis of $\PC/I(\VC(I))$:
\begin{align*}\BB=\{1,x_1,x_2,x_3,x_1 x_2,x_1 x_3,x_2 x_3\}.\end{align*}
Border basis of $I(\VC(I))$:
\begin{align*}
g_1&=1+x_2+x_3+{\bf x_1^2}\\
g_2&=-1-x_1+x_2-x_3+x_2 x_3+{\bf x_1^2 x_2}\,,\\
g_3&=-1-x_1-x_2+x_3+x_2 x_3+{\bf x_1^2 x_3}\,,\\
g_4&=3.9993-0.99984 x_1 x_2-0.99984 x_1 x_3-0.99984 x_2 x_3+{\bf x_1 x_2 x_3}\,,\\
g_5&=1+x_1+x_3+{\bf x_2^2}\,,\\
g_6&=-1+x_1-x_2-x_3+x_1 x_3+{\bf x_1 x_2^2}\,,\\
g_7&=-1-x_1-x_2+x_3+x_1 x_3+{\bf x_2^2 x_3}\,,\\
g_8&=1+x_1+x_2+{\bf x_3^2}\,,\\
g_{9}&=-1+x_1-x_2-x_3+x_1 x_2+{\bf x_1 x_3^2}\,,\\
g_{10}&=-1-x_1+x_2-x_3+x_1 x_2+{\bf x_2 x_3^2}\,.
\end{align*}
Extracted solutions (via the monomial basis):
\begin{align*}
x_1&=(-1-8.15\text{e-11}i,-1+4.37\text{e-11}i,-1-4.24\text{e-12}i)\approx (-1,-1,-1) \,,\\
x_2&=(-1.16\text{e-5}+1.41i,0.999-1.41i,0.000654+1.41i)\approx (\sqrt{2}i,1-\sqrt{2}i,\sqrt{2}i)\,,\\
x_3&=(-4.8\text{e-5}-1.41i,1+1.41i,0.000147-1.41i)\approx (-\sqrt{2}i,1+\sqrt{2}i,-\sqrt{2}i)\,,\\
x_4&=(-9.92\text{e-7}+1.41i,0.000713+1.41i,0.999-1.41i)\approx (\sqrt{2}i,\sqrt{2}i,1-\sqrt{2}i)\,,\\
x_5&=(-3.98\text{e-5}-1.41i,0.000146-1.41i,1+1.41i)\approx (-\sqrt{2}i,-\sqrt{2}i,1+\sqrt{2}i)\,,\\
x_6&=(1+1.41i,-0.000149-1.41i,-0.000145-1.41i)\approx  (1+\sqrt{2}i,-\sqrt{2}i,-\sqrt{2}i)\,,\\
x_7&=(1-1.41i,0.000103+1.41i,0.000104+1.41i)\approx  (1-\sqrt{2}i,\sqrt{2}i,\sqrt{2}i)\,.
\end{align*}

\noindent
For this example more accurate solutions:
\begin{align*}
x_1&=(-1-2.51\text{e-11}i,-1-9.85\text{e-11}i,-1-5.99\text{e-11}i)\,,\\
x_2&=(-1.73\text{e-5}+1.41i,1-1.41i,-1.29\text{e-5}+1.41i)\,,\\
x_3&=(-3.78\text{e-5}-1.41i,1+1.41i,-8.52\text{e-5}-1.41i)\,.\\
x_4&=(-2.83\text{e-5}+1.41i,4.23\text{e-5}+1.41i,1-1.41i)\,,\\
x_5&=(-4\text{e-5}-1.41i,3.88\text{e-5}-1.41i,1+1.41i)\,,\\
x_6&=(1+1.41i,-0.000191-1.41i,6.2\text{e-5}-1.41i)\,,\\
x_7&=(1-1.41i,-4.61\text{e-5}+1.41i,0.000104+1.41i)\,,
\end{align*}
could be obtain by means of the SVD-method. 
\end{example}


\section{Concluding Remarks}

In this paper we have provided a new semidefinite characterization of the real radical 
ideal of an ideal $I\subseteq\oR[x]$ as well as an algorithm to compute  all (finitely many) points
of $V_\oR(I)$ and  a set of generators (or a Gr\"obner base) of $I(V_\oR(I))$. 
The main feature of our approach is its real algebraic nature as it avoids considering complex zeros, and does not need to compute a Gr\"obner base of $I$. 
An essential step in our algorithm consists in solving the semidefinite program (\ref{P0}). Thus our algorithm is numerical.

Let us briefly mention the {\it numerical} versus {\it numeric-symbolic} (or 
{\it arbitrary precision}) issue. Some advocate that only computation with arbitrary or guaranteed precision should be
permitted while others admit some numerical imprecision; see e.g. Revol and Rouillier \cite{revol},
Stetter \cite{St04}. Clearly, being numerical in nature, the algorithm of the present paper
admits some intrinsic numerical imprecision, no matter how good are (or will be) the SDP software packages. At this stage, the only answer we propose  in this `approximate vs exact' debate is to validate or invalidate
the method by experiments. For instance,  on a significant sample, compute $J\approx I(V_\mathbb{R}(I))$ with our method
and check afterwards by symbolic methods if $\VC(J)=V_\mathbb{R}(I)$. On the other hand, the present algorithm is rather intended to illustrate that the new semidefinite characterizations of
$V_\oR(I)$ and $I(V_\oR(I))$ are
directly implementable in a relatively simple manner; clearly, its numerical features 
(like precision and stability) need further investigations beyond the scope of the present paper.

\bigskip
{\bf \large Acknowledgements.}  We are very grateful to two referees for their careful reading and their many useful suggestions that helped us improve the presentation of the paper.
In particular we thank a referee for pointing out to us that the method can also detect nonexistence of (real) solutions.
We also thank Etienne de Klerk for helpful discussions on the self-dual embedding technique for semidefinite programming and Johan L\"ofberg for his support with Yalmip.

\end{document}